\documentclass[12pt,a4paper,reqno]{amsart}
\usepackage{amsmath}
\usepackage{physics}
\numberwithin{equation}{section}
\usepackage{graphicx}
\usepackage{amssymb}
\usepackage{geometry}
\usepackage{amsthm}
\usepackage{enumitem}
\usepackage[utf8]{inputenc}
\usepackage{tikz}
\usepackage[T1]{fontenc}
\usepackage{tikz-cd}
\usepackage{txfonts}
\usepackage{mathrsfs}
\usepackage{hyperref}
\usetikzlibrary{calc}

\newtheorem{thm}{Theorem}[section]
\newtheorem*{theorem*}{Theorem}

\newtheorem{prop}[thm]{Proposition}
\newtheorem{lm}[thm]{Lemma}

\newtheorem{coro}[thm]{Corollary}

\makeatletter
\newcommand*\bigcdot{\mathpalette\bigcdot@{.5}}
\newcommand*\bigcdot@[2]{\mathbin{\vcenter{\hbox{\scalebox{#2}{$\m@th#1\bullet$}}}}}
\makeatother

\newtheorem{innercustomgeneric}{\customgenericname}
\providecommand{\customgenericname}{}
\newcommand{\newcustomtheorem}[2]{%
	\newenvironment{#1}[1]
	{%
		\renewcommand\customgenericname{#2}%
		\renewcommand\theinnercustomgeneric{##1}%
		\innercustomgeneric
	}
	{\endinnercustomgeneric}
}
\newcustomtheorem{customthm}{Theorem}
\newcustomtheorem{customcor}{Corollary}

\theoremstyle{definition}
\newtheorem{example}[thm]{Example}
\newtheorem{df}[thm]{Definition}
\newtheorem{remark}[thm]{Remark}
\newtheorem{conj}[thm]{Conjecture}

\allowdisplaybreaks

\newcommand{\N}{\mathbb{N}}
\newcommand{\Z}{\mathbb{Z}}

\newcommand{\fusion}[3]{{\binom{#3}{#1\;#2}}}

\providecommand{\keywords}[1]
{
	\small	
	\textbf{\textit{Keywords---}} #1
}

\def \ra {\rightarrow}

\def \C {\mathbb{C}}

\def\la{\lambda}

\def \al{\alpha}

\def \om{\omega}
\def \ga {\gamma}

\def \b {\beta}
\def \op {\oplus}
\def\ds{\dots}

\def \ssq{\subseteq}

\def \vac {\mathbf{1}}

\def \g {\mathfrak{g}}
\def \h {\mathfrak{h}}

\def \End {\mathrm{End}}

\def\spn{\mathrm{span}}
\def\Id{\mathrm{Id}}

\def \gr {\mathrm{gr}}
\def \wt {\mathrm{wt}}

\def \bs {\backslash}

\def\<{\langle}
\def\>{\rangle}

\title{Noetherianity of twisted Zhu algebra and bimodules}
	\author{Jianqi Liu}
\address{Department of Mathematics, University of Pennsylvania, Philadelphia, PA, 19104}
\email{jliu230@sas.upenn.edu}

\begin {document}
\maketitle

\begin{abstract}
 In this paper we show that for a large natural class of vertex operator algebras (VOAs) and their modules, the Zhu algebras and bimodules (and their $g$-twisted analogs) are Noetherian. These carry important information about the representation theory of the VOA, and its fusion rules, and the Noetherian property gives the potential for (non-commutative) algebro-geometric methods to be employed in their study.
 	\end{abstract}

 	%	\tableofcontents
 
 \section{Introduction}

This paper is concerned with finiteness property of (twisted) Zhu's algebra and its bimodules. The Zhu algebra and its twisted higher-level versions carry representation theoretical information about modules and twisted modules over VOAs. Bimodules of Zhu's algebra have been used to compute (twisted) fusion rules. It has been observed in examples that Zhu's algebras are often Noetherian and even finite-dimensional. This is unexpected given the analogy between Zhu's algebra and universal enveloping algebra of a Lie algebra since it is an open problem whether the latter is Noetherian when the Lie algebra is not finite-dimensional \cite{GW89}. By leveraging the relationships between the twisted Zhu's algebra and Zhu's $C_2$-algebra, we show the Notherianity for a large class.

In the study of the modular invariance property of vertex operator algebras (VOAs), Zhu introduced an associative algebra $A(V)$ attached to a VOA of CFT-type \cite{Z96}. Associated to an admissible $V$-module $M$, an $A(V)$-bimodule $A(M)$ was introduced by Frenkel and Zhu in order to compute the fusion rules among irreducible modules over affine VOAs \cite{FZ92}. The main result in this paper is that $A(V)$, together with its $g$-twisted analog $A_g(V)$ \cite{DLM98}, and bimodule $A_g(M)$ \cite{JJ16,GLZ23}, are left (or right) Noetherian if $V$ is $C_1$-cofinite \cite{L99,KL99} and $M$ is (weakly) $C^g_1$-cofinite. If, in addition, $V$ is $C_2$-cofinite \cite{Z96}, then the $g$-twisted higher order generalizations $A_{g,n}(V)$ \cite{DLM98(2),DLM98(3)} and $A_{g,n}(M)$ \cite{JJ16} are finite-dimensional for all $n\geq 0$. These algebraic structures encode important information about the representation theory of the VOAs including the fusion rules. Noetherianity is one of the most important finiteness properties, which gives tools for their study, for instance from (non-commutative) algebraic geometry \cite{DGK24}.

Zhu proved in \cite{Z96} that there is a one-to-one correspondence between irreducible $V$-modules and irreducible $A(V)$-modules, which leads to an equivalency between the categories of $V$-modules and $A(V)$-modules for rational VOAs. Zhu's result was generalized by Dong, Li, and Mason to the $g$-twisted case in \cite{DLM98}, and higher order (twisted) cases in \cite{DLM98(2),DLM98(3)}, wherein the notions of $g$-twisted Zhu's algebra $A_g(V)$, higher order Zhu's algebra $A_n(V)$ for $n\geq 0$, and $g$-twisted higher order Zhu's algebra $A_{g,n}(V)$ were introduced, and the one-to-one correspondences between irreducible ($g$-twisted) $V$-modules and irreducible modules over these generalized Zhu's algebra were established. From this point of view, Zhu's algebra and its generalizations tell us about the representation theory of VOAs. 

%Noetherianity is a fundamental property of a ring. It is expected to be true for any reasonable commutative ring or associative algebra. 

Dong, Li, and Mason proved that $A(V)$ is finite-dimensional if $V$ is $C_2$-cofinite \cite{DLM98,DLM00}. 
For the classical non-$C_2$-cofinite VOAs like the vacuum module VOA $V_{\hat{\g}}(\ell,0)$, the Heisenberg VOA $M_{\hat{\h}}(\ell,0)$, and the universal Virasoro VOA $\bar{V}(c,0)$ \cite{FZ92,LL04}, their Zhu's algebra are isomorphic to $U(\g)$, $\C[x_1,\ds, x_n]$, and $\C[x]$, respectively. Although these associative algebras are infinite-dimensional, they are all Noetherian. Moreover, in numerous calculations for concrete examples  \cite{FZ92,W93,DLM97,DN99,ALY14,AB23}, we see that Zhu's algebra is close to a quotient algebra of certain universal enveloping algebra of a Lie algebra $\g$. In fact, it was proved by He in \cite{He17} that the higher order Zhu's algebra $A_n(V)$ is isomorphic to a subquotient algebra of the degree zero part of the universal enveloping algebra $U(V)$ of a VOA defined by Frenkel and Zhu \cite{FZ92}. VOAs generalize Lie algebras and so the Noetherian property is unexpected given what is known about Lie algebras. For instance, if $\g$ is the Witt algebra, it was proved by Sierra and Walton that $U(\g)$ is {\em not} Noetherian \cite{SW14}. Adding to the unexpectedness of the result, $A_g(V)$ is Noetherian for all $C_1$-cofinite $V$, an unrestricted class, encompassing what are considered all reasonable examples, including the non $C_2$-cofinite VOAs mentioned above.

 The Noetherianity for Zhu's algebra has been established as an ingredient for the study of the representation theory of $C_1$-cofinite VOAs. It was used in a recent work of Damiolini, Gibney, and Krashen in \cite{DGK24}. 

To state our main results, and describe how they are proved, we introduce some notation. Let $V$ be a VOA of CFT-type, $g\in \mathrm{Aut}(V)$ be an automorphism of order $T<\infty$, and $R_2(V)=V/C_{2}(V)$ be the $C_2$-algebra \cite{Z96}. It was observed by Zhu that $A(V)$ has a filtration $\{F_pA(V)\}_{p=0}^\infty $ obtained by the grading $V=\oplus_{p=0}^\infty V_p$. The associated graded algebra $\gr A(V)$ is commutative and unital.  Arakawa,  Lam, and Yamada observed that there is an epimorphism $R_2(V)\ra \gr A(V)$ of commutative algebras \cite{ALY14}. It turns out that this epimorphism is quite useful for the study of the structure theory of $A(V)$. Using this morphism, Yang and the author proved a Schur's lemma for $C_1$-cofinite VOAs over an arbitrary field \cite{YL23}. The twisted Zhu's algebra $A_g(V)$ carries a similar level filtration $\cup_{p=0}^\infty F_p A_g(V)$, and there exists epimorphism from $R_2(V)$ to the associated graded algebra $\gr A_g(V)$ as well. It was proved by Li that $R_2(V)$ is a finitely generated algebra if $V$ is $C_1$-cofinite \cite{L99,L05}. Combining these facts together, we can prove our first main theorem (see Theorem~\ref{thm:AgVN}):
\begin{customthm}{A}\label{introthmA}
	Let $V$ be a CFT-type VOA that is $C_1$-cofinite, and let $g\in \mathrm{Aut}(V)$ be a finite order automorphism. Then $A_g(V)$ is left and right Noetherian as an associative algebra.
	\end{customthm}

The $A(V)$-bimodule $A(M)$ and its twisted analog $A_g(M)$ were introduced to compute the fusion rules among ($g$-twisted)-modules over $V$ \cite{FZ92,Li99(2),Liu23,GLZ23}. Li introduced a cofinite condition for $V$-modules, which we call the {\em wealky $C_1$-cofinite condition}, and proved that the fusion rule among three irreducible untwisted $V$-modules $M^1,M^2$, and $M^3$ is finite if $M^1$ is weakly $C_1$-cofinite, see \cite{L99}. In order to handle the $g$-twisted case, we modify Li's confinite condition and introduce a subspace $\widetilde{C}^g_1(M)$ associated to $M$. We say that $M$ is {\em weakly $C_1^g$-cofinite} if $\dim M/\widetilde{C}^g_1(M)<\infty$, see Definition~\ref{weaklyC1gcofinite}. Huang independently introduced another $C_1$-cofinite condition for modules in \cite{H05}, which is slightly stronger than Li's $C_1$-condition, to guarantee the convergence of iterated intertwining operators. As an application, Huang also proved that the fusion rule among $V$-modules $M^1,M^2$, and $M^3$ is finite if $M^1$ is $C_1$-cofinite. These $C_1$-cofinite conditions for $V$-modules correspond to finite generation properties of the twisted bimodule $A_g(M)$ over twisted Zhu's algebra $A_g(V)$. The following is our second main theorem (see Theorem~\ref{thm:AgMfg}):

\begin{customthm}{B}\label{introthmB}
Let $M$ be an untwisted admissible $V$-module. Then 
\begin{enumerate}
	\item $A_g(M)$ is finitely generated as a left or right $A_g(V)$-module if $M$ is $C_1$-cofinite.
	\item  $A_g(M)$ is finitely generated as an $A_g(V)$-bimodule if $M$ is weakly $C^g_1$-cofinite.
	\end{enumerate}
In particular, for a $C_1$-cofinite VOA $V$, $A_g(M)$ is Noetherian as a left or right $A_g(V)$-module if $M$ is $C_1$-cofnite; $A_g(M)$ is Noetherian as an $A_g(V)$-bimodule if $M$ is weakly $C^g_1$-cofinite.
\end{customthm}
 
As a Corollary of Theorem~\ref{introthmB}, using the $g$-twisted fusion rules theorem proved by Gao, the author, and Zhu in \cite{GLZ23}, we can prove following finiteness property for fusion rules among $g$-twisted modules, which simultaneously generalizes both Li and Huang's result about finiteness of fusion rules under $C_1$ condition to the $g$-twisted case (see Corollary~\ref{corofinitenessoffusion}):

\begin{customcor}{C}\label{introthmC}
Let $M^1$ be an untwisted ordinary $V$-module, and $M^2,M^3$ be $g$-twisted ordinary $V$-modules. If the $M^1$ is (weakly) $C_1$-cofinite, then the fusion rule $N\fusion{M^1}{M^2}{M^3}$ is finite. 
\end{customcor}

Theorem~\ref{introthmA} gives us a sufficient condition for the Noetherianity of $A_g(V)$. In Section~\ref{sec4}, we use a concrete example to show that the Noetherianity of $A_g(V)$ fails to be true in general if the CFT-type VOA $V$ is not $C_1$-cofinite. Since the classical examples of CFT-type VOAs (rational or not) are all $C_1$-cofinite \cite{DLM02}, it is not trivial to find a CFT-type non-$C_1$-cofinite VOA. The example we construct is a subVOA $V_M=\bigoplus_{\ga\in M} M_{\hat{\h}}(1,\ga)$ of the lattice VOA $V_{A_2}$, where $M=\{m\al+n\b: m\geq n\geq 1\}\cup \{0\}$ is an abelian submonoid of the root lattice $A_2$, see Figure~\ref{fig1}. This example is a modification of the Borel-type subVOA of a lattice VOA defined by the author in \cite{Liu24}. Using a similar method as in \cite{Liu24}, we can give an explicit description of the Zhu's algebra $A(V_M)$ of $V_M$ and show that it is not Noetherian. The following is our third main theorem (see Theorem~\ref{thm:V_MisnotC1}, Theorem~\ref{thm4.6}, and Corollary~\ref{corononnoehterian}): 

\begin{customthm}{D}\label{introthmD} Let $V_M=M_{\hat{\h}}(1,0)\op\bigoplus_{m\geq n\geq 1} M_{\hat{\h}}(1,m\al+n\b)$. Then 
\begin{enumerate}
	\item $V_M/C_1(V_M)$ has a basis $\{\vac+C_1(V_M),e^{m\al+\b}+C_1(V_M):m\geq 1\}$. In particular, the CFT-type VOA $V_M$ is {\em  not} $C_1$-cofinite. 
	\item  $A(V_M)\cong \C[x,y]\op \left(\bigoplus_{m=1}^\infty z_m\C[y]\right),$ where $J=\bigoplus_{m=1}^\infty z_m\C[y]$ is a two-sided ideal of $A(V_M)$ which is not finitely generated. In particular, $A(V_M)$ is not Noetherian. 
\end{enumerate}	
	 
	\end{customthm}

The higher level generalization of Zhu's algebra $A_n(V)$ was introduced by Dong, Li, and Mason in \cite{DLM98(3)} to study the rationality of VOAs. They proved that $V$ is rational if and only if $A_n(V)$ are semisimple for all $n\geq 0$. $A_n(V)$ was further generalized to the $g$-twisted case in \cite{DLM98(2)}. The $g$-twisted higher Zhu's algebra $A_{g,n}(V)$ controls the first $n$ level of a $g$-twisted admissible $V$-module $M$, where $n=l+\frac{i}{T}$ with $l\in \N$ and $0\leq i\leq T-1$. In Section~\ref{sec4.1}, we introduce a shifted level-filtration $\cup_{p=2l}^\infty F_p A_{g,n}(V)$ on $A_{g,n}(V)$ which is compatible with the product on $A_{g,n}(V)$, see Lemma~\ref{lmgrAgnV}. On the other hand, Zhu's $C_2$-algebra $R_2(V)$ also has a higher order generalization $R_{2l+2}(V)=V/C_{2l+2}(V)$. However, unlike $R_2(V)$, the  associative algebra $R_{2l+2}(V)$ is not commutative in general. In Section~\ref{sec4.2},  we show that there is a surjective linear map $\varphi_n: R_{2l+2}(V)\ra \gr A_{g,n}(V)$, which is a homomorphism of associative algebras if $i<\lfloor T/2\rfloor$, see Theorem~\ref{mainthm2}

Gaberdiel and Neitzke proved that the $C_2$-cofinite condition is strong enough so that it implies $\dim R_{2l+2}(V)<\infty$ for all $l\geq 0$, see \cite{GN03}. Using this fact, Miyamoto proved that $A_n(V)$ are finite-dimensional for all $n\geq 0$ if $V$ is $C_2$-cofinite, which is a key property for the modular invariance of pseudo trace functions of $C_2$-cofinite VOAs \cite{M04}. Buhl found a module version of Gaberdiel and Neitzke's theorem and proved that $A_n(M)$ are finite-dimensional for all $n\geq 0$ if $V$ is $C_2$-cofinite and $M$ is $C_2$-cofinite \cite{Bu02}. The finiteness of $A_n(M)$ could be useful in generalizing Huang's modular invariance of logarithmic intertwining operators \cite{H23} to $C_2$-cofinite but not necessarily rational VOAs. 
With the surjective linear map $\varphi_n: R_{2l+2}(V)\ra \gr A_{g,n}(V)$, we can prove our last main theorem, which is a twisted version of Miyamoto and Buhl's finiteness results about $A_n(V)$ and $A_n(M)$ (see Corollary~\ref{coroofmain2} and Corollary~\ref{coro2ofmain2}): 

\begin{customthm}{E}\label{introthmE}
	 Let $M$ be an untwisted irreducible admissible $V$-module, and let $n=l+\frac{i}{T}\in \frac{1}{T}\Z$, where $l\in \N$ and $0\leq i\leq T-1$.
	\begin{enumerate}
		\item If $V$ is $C_2$-cofinite, then $ A_{g,n}(V)$ is a finite-dimensional associative algebra, and $A_{g,n}(M)$ is a finite-dimensional $A_{g,n}(V)$-bimodule. 
		\item If $i< \lfloor T/2\rfloor$, and $R_{2l+2}(V)$ is a finitely generated associative algebra, then $ A_{g,n}(V)$ is left and right Noetherian. If, furthermore, $M$ is $C_{2l+1}$-cofinite, then $ A_{g,n}(M)$ is left and right Noetherian.
	\end{enumerate}
	\end{customthm}

We conjecture that $A_{g,n}(V)$ are left and right Noetherian for all $n\geq 0$ if $V$ is $C_1$-cofinite. According to a recent structural result about the higher order Zhu's algebra of the Heisenberg VOA by Damiolini, Gibney, and Krashen in \cite{DGK23}, we know that this conjecture is true if $V$ is the Heisenberg VOA and $g=\Id_V$.

This paper is organized as follows: In Section~\ref{Sec2} we recall the definitions of $g$-twisted modules, twisted Zhu's algebra $A_g(V)$ and its bimodule $A_g(M)$, the $C_2$-algebra $R_2(V)$ and its relation with the $C_1$-cofinite condition. In Section~\ref{Sec3}, we prove Theorem~\ref{introthmA}, Theorem~\ref{introthmB}, and Corollary~\ref{introthmC}. In Section~\ref{sec4}, we introduce the CFT-type VOA $V_M$ and prove that it is not $C_1$-cofinite. Then we determine the Zhu's algebra $A(V_M)$ and show that it is not Noetherian as claimed in Theorem~\ref{introthmD}. 
 In Section~\ref{Sec4}, we first introduce a shifted level filtration on $A_{g,n}(V)$ and discuss its relations with the $C_{2l+2}$-algebra $R_{2l+2}(V)$, then we use it to prove Theorem~\ref{introthmE}.

{\bf Convention}: All vector spaces are defined over $\C$, the field of complex number. $\N$ represents the set of natural numbers including $0$.  
 
 \keywords{{\bf Keywords}: vertex operator algebra, associative algebra, twisted representation, fusion rules, Noether ring}

\section{Preliminaries}\label{Sec2}

\subsection{$g$-twisted modules over vertex operator algebras}
For the definitions of vertex operator algebras (VOAs), untwisted modules over VOAs, Zhu's algebra and its bimodule, we refer to \cite{FLM88,FZ92,FHL93,DL93,Z96,LL04}. Throughout this paper, we assume a VOA $(V,Y,\vac,\om)$ is of {\em CFT-type}: $V=V_{0}\oplus V_{+}$, where $V_{0}=\C\vac$ and $V_{+}=\bigoplus_{n=1}^{\infty}V_{n}$. 

Let $g:V\ra V$ be an automorphism of $V$ finite order $T$ \cite{FLM88}. Then $V$ has a $g$-eigenspace decomposition \cite{DL93,DVVV89}: 
\begin{equation}\label{Veigendec}
V=\bigoplus_{r=0}^{T-1}V^r,\quad \mathrm{where}\quad  V^r=\{a\in V: g(a)=e^{\frac{2\pi i r}{T}}a\}.
\end{equation}
In the rest of this paper, we fix an automorphism $g\in \mathrm{Aut}(V)$ of order $T$.
\begin{df}{\cite{DL93,DLM98,H10}}\label{def:twistedmodule}\label{def:gtwisted}
	 A {\em $g$-twisted weak $V$-module} is a pair $(M,Y_M)$, where $M$ is a vector space, and $Y_M$ a linear map 
	\begin{align*}
	Y_M:V&\ra \End(M)[[z^{1/T},z^{-1/T}]]\\
	a&\mapsto Y_M(a,z)=\sum_{n\in \Z} a_n z^{-n-1-\frac{r}{T}},\quad \mathrm{for}\ a\in V^r,
	\end{align*}
	satisfying the following properties:
\begin{enumerate}[label=(\alph*)]
	\item {\em(truncation property)} For any $a\in V$ and $v\in M$, we have $a_nv=0$ for $n\in \frac{1}{T}\Z$ and $n\gg 0$. 
	\item {\em ($g$-twisted Jacobi identity)} For any $a\in V^r$ with $0\leq r\leq T-1$, and $b\in V$, we have 
	\begin{equation}
	\begin{aligned}
	&z_{0}^{-1}\delta\left(\frac{z_1-z_2}{z_0}\right) Y_M(a,z_1)Y_M(b,z_2)-	z_{0}^{-1}\delta\left(\frac{-z_2+z_1}{z_0}\right)Y_M(b,z_2)Y_M(a,z_1)\\
	&=	z_{2}^{-1}\delta\left(\frac{z_1-z_0}{z_2}\right)\left(\frac{z_1-z_0}{z_2}\right)^{-r/T} Y_M(Y(a,z_0)b,z_2). 
	\end{aligned}
	\end{equation}
	\end{enumerate}
A $g$-twisted weak $V$-module $M$ is called {\em admissible} if $M$ has a subspace decomposition: $$M=\bigoplus_{n\in \frac{1}{T}\N}M(n),$$ such that $a_mM(n)\ssq M(\wt a-m-1+n)$ for all $a\in V$ homogeneous, $m\in \frac{1}{T}\Z$, and $n\in \frac{1}{T}\N$. 

An admissible $g$-twisted $V$-module $M$ is called an {\em (ordinary) $g$-twisted $V$-module} if there exists $\la \in  \C$, called the {\em conformal weight}, such that $M(n)=M_{\la+n}$ is an eigenspace of $L(0)$ of eigenvalue $\la+n$, and $M(n)$ is finite-dimensional, for all $n\in \frac{1}{T}\N$. 
	\end{df}

In particular, if $g=\Id_V$ and $T=1$, then Definition~\ref{def:twistedmodule} recovers the usual definitions of weak $V$-modules, admissible $V$-modules, and ordinary $V$-modules. 
 
 \subsection{The $g$-twisted Zhu's algebra $A_{g}(V)$ and it bimodule $A_{g}(M)$}

The $g$-twisted Zhu's algebra $A_g(V)$ was constructed by Dong, Li, and Mason in \cite{DLM98}, as a $g$-twisted generalization of the usual Zhu's algebra $A(V)$ in \cite{Z96}, which controls the bottom level $M(0)$ of a $g$-twisted admissible $V$-module. 

\subsubsection{Definition of $A_g(V)$} By definition, for any $a\in V^r$ with $0\leq r\leq T-1$, and $b\in V$, let
\begin{equation}\label{def:OgV}
a\circ_{g} b:= \Res_z Y(a,z)b\frac{(1+z)^{\wt a-1+\delta(r)+\frac{r}{T}}}{z^{1+\delta(r)}},\quad \mathrm{where}\quad \delta(r)=\begin{cases}1&\mathrm{if}\ r=0\\ 0&\mathrm{if}\ r>0 \end{cases}.
\end{equation}

Let $O_g(V):=\spn\{a\circ_g b: a\in V^r,\ 0\leq r\leq T-1,\ b\in V \}$, and $A_g(V):=V/O_g(V)$. Define 

\begin{equation}\label{def:stargV}
a\ast_g b:=\begin{cases}\Res_z Y(a,z)b\frac{(1+z)^{\wt a}}{z}=\sum_{j\geq 0}\binom{\wt a}{j} a_{j-1}b &\mathrm{if}\ a\in V^0\\ 0&\mathrm{if}\ a\in V^r,\ r>0. \end{cases}
\end{equation}
By Theorem 2.4 in \cite{DLM98}, $A_g(V)$ is an associative algebra with respect to the product \eqref{def:stargV}, with unit element $[\vac]=\vac+O_g(V)$. By Lemma 2.2 in \cite{DLM98}, we have 
\begin{equation}\label{commprodAgV}
a\ast_g b-b\ast_g a=\begin{cases}\Res_z Y(a,z)b(1+z)^{\wt a-1}=\sum_{j\geq 0}\binom{\wt a-1}{j} a_jb &\mathrm{if}\ a\in V^0\\ 0&\mathrm{if}\ a\in V^r,\ r>0. \end{cases}
\end{equation}

\subsubsection{Definition of $A_g(M)$} Let $M$ be a (untwisted) admissible $V$-module. The $A_g(V)$-bimodule $A_g(M)$ was first introduced in \cite{JJ16} as a $g$-twisted generalization of the $A(V)$-bimodule $A(M)$ in \cite{FZ92}. One can use $A_g(M)$ and $A_g(V)$ to calculate the fusion rules among one untwisted $V$-module $M^1$ and two $g$-twisted $V$-modules $M^2$ and $M^3$, see \cite{GLZ23}. 

Similar to \eqref{def:OgV}, for any $a\in V^r$ with $0\leq r\leq T-1$, and $v\in M$, we let 
\begin{equation}\label{def:OgM}
a\circ_g v:=\Res_z Y_M(a,z)v\frac{(1+z)^{\wt a-1+\delta(r)+\frac{r}{T}}}{z^{1+\delta(r)}}. 
\end{equation}
Let $O_g(M)=\spn\{a\circ_gv: a\in V^r,\ 0\leq r\leq T-1,\ v\in M \}$, and $A_g(M)=M/O_g(M)$. Define:

\begin{align}
a\ast_g v:&=\begin{cases}\Res_z Y_M(a,z)v\frac{(1+z)^{\wt a}}{z}=\sum_{j\geq 0}\binom{\wt a}{j} a_{j-1}v &\mathrm{if}\ a\in V^0\\ 0&\mathrm{if}\ a\in V^r,\ r>0, \end{cases}\label{leftmodact}\\
v\ast_g a:&=\begin{cases}\Res_z Y_M(a,z)v\frac{(1+z)^{\wt a-1}}{z}=\sum_{j\geq 0}\binom{\wt a-1}{j} a_{j-1}v &\mathrm{if}\ a\in V^0\\ 0&\mathrm{if}\ a\in V^r,\ r>0. \end{cases}\label{rightmodact}
\end{align}

Then $A_g(M)$ is a bimodule over $A_g(V)$ with respect to the left and right actions \eqref{leftmodact} and \eqref{rightmodact}, see \cite{JJ16} Theorem 3.4 or \cite{GLZ23} Lemma 6.1. The following formula follows immediately from \eqref{leftmodact} and \eqref{rightmodact}:  

\begin{equation}\label{commAgM}
a\ast_g v-v\ast_g a=\begin{cases}\Res_z Y_M(a,z)v (1+z)^{\wt a-1}=\sum_{j\geq 0}\binom{\wt -1}{j} a_jv &\mathrm{if}\ a\in V^0\\ 0&\mathrm{if}\ a\in V^r,\ r>0. \end{cases}
\end{equation}
Moreover, using the $L(-1)$-derivative property of $Y_M$, one can show 
\begin{equation}\label{propertyofOgM}
\Res_z Y_M(a,z)v\frac{(1+z)^{\wt a-1+\delta(r)+\frac{r}{T}+n}}{z^{1+\delta(r)+m}}\in O_g(M),\quad m\geq n\geq 0. 
\end{equation}

\subsubsection{Level filtration on $A_g(V)$ and $A_g(M)$} For the general theory of filtered rings and modules, we refer to \cite{MR87}. It was observed by Zhu in \cite{Z96} that $A(V)$ has a canonical filtration obtained by the level decomposition of $V$: 
$$A(V)=\bigcup_{p=0}^\infty F_pA(V),\quad \mathrm{where}\quad  F_pA(V)=\left(\oplus_{n=0}^p V_n+O(V)\right)/O(V).$$
 
We can similarly define the level filtration on $A_g(V)$ and $A_g(M)$ as follows: 
\begin{align}
A_g(V)&=\bigcup_{p=0}^\infty F_pA_g(V),\quad \mathrm{where}\quad  F_pA_g(V)=\left(\oplus_{n=0}^p V_n+O_g(V)\right)/O_g(V).\label{levelfilAgV}\\
A_g(M)&=\bigcup_{p=0}^\infty F_pA(M),\quad \mathrm{where}\quad  F_pA_g(M)=\left(\oplus_{n=0}^p M(n)+O_g(M)\right)/O_g(M).\label{levelfilAgM}
\end{align}
 
\begin{lm}\label{lm:grAgV}
	Let $V$ be a VOA, and $M$ be an admissible untwisted $V$-module. Then 
	\begin{enumerate}
\item  $A_g(V)$ is a filtered associated algebra with respect to the filtration \eqref{levelfilAgV}, and the associated graded algebra 
$$\gr A_g(V)=\bigoplus_{p=0}^\infty F_{p}A_g(V)/F_{p-1}A_g(V)\quad \mathrm{with}\quad F_{-1}A_g(V)=0$$ is a commutative associative algebra with respect to the product: 
\begin{equation}\label{prodongrAV}
\left([a]+F_{p-1}A_g(V)\right)\ast_g \left([b]+F_{q-1}A_g(V)\right)=\begin{cases}[a_{-1}b]+F_{p+q-1}A_g(V) &\mathrm{if}\ a\in V^0,\\
0+F_{p+q-1}A_g(V)&\mathrm{if}\ a\in V^r,\ r>0,
\end{cases}
\end{equation}
for any $a\in \oplus_{n=0}^p V_n$ and $b\in \oplus_{n=0}^q V_n$, and $p,q\geq 0$, with identity element $[\vac]\in F_{0}A_g(V)$.
\item 	 $A_g(M)$ is a filtered $A_g(V)$-bimodule with respect to \eqref{levelfilAgV} and \eqref{levelfilAgM}, and the associated graded space 
$$\gr A_g(M)=\bigoplus_{p=0}^\infty F_pA_g(M)/F_{p-1}A_g(M)\quad \mathrm{with}\quad F_{-1}A_g(M)=0$$
 is a graded $\gr A_g(V)$-module with respect to the module action:
 \begin{equation}\label{prodongrAM}
 \left([a]+F_{p-1}A_g(V)\right)\ast_g \left([v]+F_{q-1}A_g(M)\right)=\begin{cases}[a_{-1}v]+F_{p+q-1}A_g(M) &\mathrm{if}\ a\in V^0,\\
 0+F_{p+q-1}A_g(M)&\mathrm{if}\ a\in V^r,\ r>0,
 \end{cases}
 \end{equation}
 for any $a\in \oplus_{n=0}^p V_n$ and $v\in \oplus_{n=0}^q M(n)$, and $p,q\geq 0$. 
\end{enumerate}
	\end{lm}
\begin{proof}
 By \eqref{def:stargV}, it is clear that $F_pA_g(V)\ast_g F_qA_g(V)\ssq F_{p+q}A_g(V)$ for any $p,q\geq 0$, since we have $[a]\ast_g[b]=\sum_{j\geq 0}\binom{\wt a}{j} [a_{j-1} b]$ or $0$, and $\wt (a_{j-1}b)=\wt a-j+\wt b\leq p+q$ for $a\in \oplus_{n=0}^p V_n$ and $b\in \oplus_{n=0}^q V_n$. Thus, $A_g(V)$ is a filtered algebra, and $\gr A_g(V)$ is a graded algebra with respect to the product \eqref{prodongrAV}. Assume $a\in V^0$. By \eqref{commprodAgV} we have 
 \begin{align*}
 &\left([a]+F_{p-1}A_g(V)\right)\ast_g \left([b]+F_{q-1}A_g(V)\right)-\left([b]+F_{p-1}A_g(V)\right)\ast_g \left([a]+F_{q-1}A_g(V)\right)\\
 &=\sum_{j\geq 0}[a_{j}b]+F_{p+q-1}A_g(V)=0,
 \end{align*}
 since $\wt a_{j}b=\wt a-j-1+\wt b<p+q$ and so $[a_jb]\in F_{p+q-1}A_g(V)$ for all $j\geq 0$. If $a\in V^r$ with $r>0$, clearly $[a]+F_{p-1}A_g(V)$ commutes with any other elements in $\gr A_g(V)$. Thus, $\gr A_g(V)$ is a commutative associative algebra.  
 
 Similarly, by \eqref{leftmodact} and \eqref{rightmodact}, we have $F_pA_g(V)\ast_g F_qA_g(M)\ssq F_{p+q}A_g(M)$ and $F_pA_g(M)\ast_g F_qA_g(V)\ssq F_{p+q}A_g(M)$. Thus, $A_g(M)$ is a filtered $A_g(V)$-bimodule, and $\gr A_g(M)$ is a $\gr A_g(V)$-bimodule. By \eqref{commAgM}, the left and right $\gr A_g(V)$-module actions on $\gr A_g(M)$ coincide. Hence $\gr A_g(M)$ is a graded $\gr A_g(V)$-module with respect to \eqref{prodongrAM}. 
	\end{proof}
 
\subsection{The cofinite conditions of a VOA} 
The $C_2$-cofinite condition of $V$ was introduced by Zhu in \cite{Z96} to guarantee the convergence of the $n$-point trace functions. By definition, $C_2(V):=\spn\{ a_{-2}b: a,b\in V\}$, and $V$ is called $C_2$-{\em cofinite} if $\dim V/C_2(V)<\infty$. Zhu also proved in \cite{Z96} that 
$$R_2(V)=V/C_2(V)=\bigoplus_{p=0}^\infty V_p/(C_2(V)\cap V_p)$$
is a unital graded commutative associative algebra with respect to the product 
\begin{equation}\label{RVprod}
(a+C_2(V))\cdot (b+C_2(V))=a_{-1}b+C_2(V),\quad a,b\in V,
\end{equation}
with identity element $\vac+C_2(V)$.

 The notion of a strongly generated vertex operator algebra was introduced by Kac \cite{K97}:
 \begin{df}\label{df2.4}
 	Let $V$ be a VOA, and $U\ssq V$ be a subset. $V$ is said to be {\em strongly generated by $U$} if $V$ is spanned by elements of the form:
 	$$a^{1}_{-n_{1}}\ds a^{r}_{-n_{r}}u,$$
 	where $a^{1},\ds, a^{r},u\in U$, and $n_{i}\geq 1$ for all $i$. If $V$ is strongly generated by a finite dimensional subspace, then $V$ is called {\em strongly finitely generated}.
 \end{df}
 In the study of the strong generation property and PBW-basis of VOAs, Li introduced a similar condition in \cite{L99,KL99}, called $C_1$-cofiniteness. By definition, 
 \begin{equation}\label{def:C_1}
 C_1(V)=\spn\{a_{-1}b: a, b\in V_+\}+\spn\{L(-1)c: c\in V\},
 \end{equation}
 and $V$ is called {\em $C_1$-cofinite} if $\dim V/C_1(V)<\infty$.
 It is clear that $C_2(V)\ssq C_1(V)$. Hence any $C_2$-cofinite VOA is also $C_1$-cofinite.
The following Theorem that relates the $C_1$-cofinite condition with the strong generation property of a VOA was proved by Li, see \cite{L99,KL99,L05}:

 \begin{thm}\label{equiC1}
 	Let $V$ be a VOA, and $U\subseteq V_{+}$ be a graded subspace. The following conditions are equivalent:
 	\begin{enumerate}
 		\item $V$ is strongly generated by $U$. 
 		\item $V_{+}=U+C_{1}(V)$ as vector spaces.
 		\item $(U+C_{2}(V))/C_{2}(V)$ generates $R_2(V)$ as commutative algebra.
 	\end{enumerate}
 In particular, $V$ is strongly finitely generated if and only if $V$ is $C_1$-cofinite, if and only if $R_2(V)$ is a finitely generated commutative algebra. 
 \end{thm}

The $C_1$-cofiniteness condition for $V$-modules was introduced by Huang in \cite{H05}. By definition, an admissible $V$-module $M$ is called $C_1$-cofinite if $\dim M/C_1(M)<\infty$, where \begin{equation}\label{C1M}
C_1(M):=\spn\{ a_{-1}v: a\in V_+,v\in M\}.
\end{equation}
There is a similar subspace of $M$ introduced by Li in \cite{L99}: 
\begin{equation}\label{weakC1M}
B(M):=\spn\{ a_{-1}v: a\in V_+,v\in M\}+\spn\{b_0u:b\in \oplus_{n\geq 2}V_n, u\in M \}.
\end{equation}

We need to adjust the definition of $B(M)$ a little bit to make it compatible with $A_g(M)$: 
\begin{equation}\label{weakC1gM}
\widetilde{C}^g_1(M):=\spn\{ a_{-1}v: a\in V_+,v\in M\}+\spn\{b_0u:b\in \oplus_{n\geq 2}V_n\cap V^0, u\in M \},
\end{equation}
where $V^0\subset V$ is the fixed point subVOA \eqref{Veigendec} with respect to $g$. Then $\widetilde{C}_1^g(M)=B(M)$ if $g=\Id_V$. Observe that the space spanned by $b_0u$ in \eqref{weakC1gM} is nonzero, since $\om \in \oplus_{n\geq 2}V_n\cap V^0$ and there exists $u\in M\bs\{0\}$ such that $\om_0u=L(-1)u\neq 0$, see \cite{Li94}. 
\begin{df}\label{weaklyC1gcofinite}
Let $M$ be an admissible untwisted $V$-module. We say that $M$ is {\em weakly $C_1$-cofinite} if $\dim M/B(M)<\infty$; we say that $M$ is
 {\em weakly $C^g_1$-cofinite} if $\dim M/ \widetilde{C}^g_1(M)<\infty$.  
	\end{df}
Since $C_1(M)\ssq \widetilde{C}^g_1(M)\ssq B(M)$, the following Lemma is evident:

\begin{lm}
	Let $M$ be an admissible $V$-module. If $M$ is $C_1$-cofinite, then it must be weakly $C_1^g$-cofinite. If $M$ is weakly $C_1^g$-cofinite, then it must be weakly $C_1$-cofinite. 
	\end{lm}

 However, the converse of these statements is not true. 
\begin{example}
	Let $V=\bar{V}(c,0)$ be the universal Virasoro VOA with central charge $c>0$. It is well-known that the Verma module $M=M(c,h)$ over the Virasoro Lie algebra of central charge $c$ and highest weight $h>0$ is an admissible module over $V$, see \cite{FZ92,LL04}. Recall that
	$$M(c,h)=\spn\{L(-n_1)L(-n_2)\ds L(-n_n)v_{c,h}: n_1\geq n_2\geq \ds\geq n_k\geq 1\}.$$
	Then $M(c,h)=\C v_{c,h}+B(M(c,h))$ in view of \eqref{weakC1M}. Hence $M(c,h)$ is weakly $C_1$-cofinite. However, $M(c,h)=\spn\{L(-1)^kv_{c,h}:k\geq 0\}+C_1(M(c,h))$, and  $L(-1)^kv_{c,h}\neq 0$ for all $k\geq 0$ in a Verma module. Thus, $M(c,h)$ is {\em not} $C_1$-cofinite. 
	\end{example}

Finally, we recall the following fact about the Noetherianity of a filtered ring, see Theorem 6.9 in \cite{MR87}. 
 \begin{prop}\label{prop2.2}
 	Let $R$ be a filtered ring such that the associated graded ring $\mathrm{gr}R$ is left (resp. right) Noetherian, then $R$ is left (resp. right) Noetherian. 
 \end{prop}

\section{Noetherianity of twisted Zhu's algebra and its bimodule}\label{Sec3}

We prove our main theorem of this paper in this Section.

\subsection{Noetherianity of $A_g(V)$ for $C_1$-cofinite VOA $V$}
In the study of Zhu's algebra of the parafermion VOAs \cite{ALY14}, Arakawa,  Lam, and Yamada introduced an epimorphism of commutative associative algebras: 
\begin{equation}\label{gradedmap}
\begin{aligned}
\phi: R_2(V)&\ra \gr A(V)=\bigoplus_{p=0}^{\infty}F_pA(V)/F_{p-1}A(V),\\
a+C_{2}(V)&\mapsto [a]+F_{p-1}A(V),\quad   a\in \oplus_{n=0}^p V_n,
\end{aligned}
\end{equation}
where $\gr A(V)$ is the graded algebra in \ref{lm:grAgV} with $g=\Id_V$. This map was also used to prove the Schur's Lemma for irreducible modules of VOAs over an arbitrary field in \cite{YL23}. %Using a similar epimorphism for the $g$-twisted Zhu's algebra, we can prove the following theorem: 
%Define $\phi: V=\oplus_{n=0}^{\infty}V_{n}\ra \gr A(V): \ \phi(x_{1}+...+x_{r})=\overline{x_{1}}+...\overline{x_{r}}$, where $x_{i}\in V_{n_{i}}$ and $\overline{x_{i}}\in A(V)_{n_{i}}/A(V)_{n_{i}-1}$ for all $i$. Clearly $\phi$ is linear, and we claim that $\phi(C_{2}(V))=0$. \par 

\begin{thm}\label{thm:AgVN}
Let $V$ be a VOA, and $g\in \mathrm{Aut}(V)$ be a finite order automorphism. If $V$ is strongly finitely generated, or equivalently, $C_{1}$-cofinite, then $A_g(V)$ is left and right Noetherian. 
	\end{thm}
\begin{proof}
First, we generalize the epimorphism \eqref{gradedmap} to the following $g$-twisted case:
\begin{equation}\label{varphi}
\begin{aligned}
\varphi: R_2(V)&\ra \gr A_g(V)=\bigoplus_{p=0}^{\infty}F_pA_g(V)/F_{p-1}A_g(V),\\
a+C_{2}(V)&\mapsto [a]+F_{p-1}A_g(V),\quad   a\in \oplus_{n=0}^p V_n.
\end{aligned}
\end{equation}
 For any $a\in V_p$ and $b\in V_q$ with $p,q\geq 0$, we have $[a_{-2}b]\in F_{p+q+1}A_g(V)$. To show $\varphi$ is well-defined,
we need to show $[a_{-2}b]\equiv 0\pmod{F_{p+q}A_g(V)}$. We may also assume $a\in V^r$ for some $0\leq r\leq T-1$. If $r=0$, by \eqref{def:OgV} we have 
$[a\circ_g b]=\sum_{j\geq 0}\binom{\wt a}{j} [a_{j-2}b]=[0]$ in $ A_g(V).$ Hence $$[a_{-2}b]=-\sum_{j\geq 1}\binom{\wt a}{j} [a_{j-2}b]\in F_{p+q}A_g(V)$$
since $\wt (a_{j-2}b)=p-j+1+q\leq p+q$. If $r>0$,  by Lemma 2.2 in \cite{DLM98}, we have
$$\Res_z Y(a,z)b\frac{(1+z)^{\wt a-1+\frac{r}{T}}}{z^2}=\sum_{j\geq 0}\binom{\wt a-1+\frac{r}{T}}{j} a_{j-2}b\in O_g(V).$$
Hence $[a_{-2}b]=-\sum_{j\geq 1}\binom{\wt a-1+\frac{r}{T}}{j} [a_{j-2}b]\in F_{p+q}A_g(V)$, and so $\varphi$ is well-defined. 

Clearly, $\varphi$ is surjective and grading-preserving. Next, we show that $\varphi$ is a homomorphism of commutative algebras. Let $a\in V^0\cap V_p$ and $b\in V_q$, by \eqref{RVprod} and  \eqref{prodongrAV}, we have
\begin{align*}
\varphi((a+C_2(V))\cdot (b+C_2(V)))&=\varphi(a_{-1}b+C_2(V))=[a_{-1}b]+F_{p+q-1}A_g(V)\\
&=\left([a]+F_{p-1}A_g(V)\right)\ast_g\left([b]+F_{q-1}A_g(V)\right)\\
&=\varphi(a+C_2(V))\ast_g \varphi(b+C_2(V)).  
\end{align*}
Now let $a\in V^r\cap V_p$ and $b\in V_q$, for some $1\leq r\leq T-1$. By \eqref{def:OgV} we have 
$$a_{-1}b\equiv -\sum_{j\geq 1}\binom{\wt a-1+\frac{r}{T}}{j} a_{j-1}b\pmod{O_g(V)}.$$
Then $[a_{-1}b]+F_{p+q-1}A_g(V)=-\sum_{j\geq 1}\binom{\wt a-1+\frac{r}{T}}{j}[a_{j-1}b]+F_{p+q-1}A_g(V)=0+F_{p+q-1}A_g(V)$ since $\wt (a_{j-1}b)=p+q-j\leq p+q-1$ for any $j\geq 1$. Thus,
\begin{align*}
\varphi(a+C_2(V)\cdot (b+C_2(V)))=[a_{-1}b]+F_{p+q-1}A_g(V)=0=\varphi(a+C_2(V))\ast_g \varphi(b+C_2(V)),
\end{align*}
in view of  \eqref{prodongrAV}. Hence $\varphi$ in \eqref{varphi} is an grading-preserving epimorphism of commutative algebras. Since $V$ is strongly finitely generated, there exists a subspace $U=\spn\{a^{1},\ds ,a^{m}\}\subset V$ of homogeneous elements $a^1\in V_{p_1},\ds, a^m\in V_{p_m}$ that strongly generates $V$. By Theorem~\ref{equiC1}, $R(V)$ is generated by $\{a^{1}+C_{2}(V),\ds, a^{m}+C_{2}(V)\}$ as a commutative algebra. Since $\varphi$ is an epimorphism, $\gr A_g(V)$ is generated by $\{[a^1]+F_{p_1-1}A_g(V),\ds , [a^m]+F_{p_m-1}A_g(V) \}$ as a commutative algebra. 
In particular, $\gr A_g(V)$ is Noetherian since it is quotient ring of the polynomnial ring $\C[T_1,\ds, T_m]$. Then by Proposition\ref{prop2.2}, $A_g(V)$ is also left Noetherian. 
	\end{proof}

\begin{remark}
If $g=\Id_V$, we have $A_g(V)=A(V)$. Note that the conclusion in Theorem~\ref{thm:AgVN} does not depend on the choice of $g$. Thus, $A(V)$ is Noetherian if $V$ is $C_1$-cofinite. 
	\end{remark}

\subsection{ Noetherianity of $A_g(M)$ for (weakly) $C_1$-cofinite $V$-module $M$}
Let $M$ be an admissible untwisted $V$-module. It was proved by Li that if $M=W+\widetilde{C}_1(M)$, then $A(M)$ is generated by $(W+O(M))/O(M)$ as an $A(V)$-bimodule, see  \cite{L99} Proposition 3.16. We have a similar result about $A_g(M)$, combined with Huang's $C_1$-cofinite condition \eqref{C1M}. 

\begin{thm}\label{thm:AgMfg}
	Let $M$ be an admissible untwisted $V$-module. 
	\begin{enumerate}
		\item Let $M=U+C_1(M)$ and $U=\spn\{u^i:i\in I\}$. Then 
		\begin{equation}\label{AgMgenerate}
		A_g(M)=\sum_{i\in I} A_g(V)\ast_g [u^i]=\sum_{i\in I} [u^i]\ast_g A_g(V)
		\end{equation}
		as a left or right $A_g(V)$-module. In particular, $A_g(M)$ is Noetherian as a left or right $A_g(V)$-module if $V$ is $C_1$-cofinite and $M$ is $C_1$-cofinite. 
		\item Let $M=W+\widetilde{C}^g_1(M)$ and $W=\spn\{w^j:j\in J\}$. Then 
		\begin{equation}\label {AgMbigenerate}
		A_g(M)=\sum_{j\in J} A_g(V)\ast_g [w^j]\ast_g A_g(V)
		\end{equation}
		as an $A_g(V)$-bimodule.  In particular, $A_g(M)$ is Noetherian as an $A_g(V)$-bimodule if $V$ is $C_1$-cofinite and $M$ is weakly $C^g_1$-cofinite. 
		\end{enumerate}
	\end{thm}
\begin{proof}
(1) Denote the right submodule $\sum_{i\in I} [u^i]\ast_g A_g(V)$ of $A_g(M)$ by $N$. We use induction on degree $n$ of $M(n)$ to show $[M(n)]\ssq N$ in $A_g(M)$. Since $\deg (a_{-1}v)=\wt a+\deg v\geq 1$ for any $a\in V_+$, we have $C_1(M)\ssq \oplus_{m\geq 1}M(m)$. So $M(0)\ssq U$ and $[M(0)]\ssq N$. Suppose the conclusion holds for smaller $n$. Let $x\in M(n)$. We may assume
$$x=u+\sum_{k=1}^s a^k_{-1}v^k,\quad u\in U,\ a^k\in V_+\cap V^r,\ 0\leq r\leq T-1,\ v^k\in M,$$
with $\wt a^k+\deg  v^k=n$ for all $k$. Since $[u]\in N$, we need to show $[a^k_{-1}v^k]\in N$ for all $1\leq k\leq s$. 

Fix a $1\leq k\leq s$. If $r=0$, by \eqref{rightmodact} we have 
\begin{equation}\label{eqmiddlestep1}
[a^k_{-1}v^k]=[v^k]\ast_g [a^k]-\sum_{j\geq 1}\binom{\wt a-1}{j}[a^k_{j-1}v^k].
\end{equation}
Note that $\deg v^k<n$ since $\wt a^k\geq 1$. By the induction hypothesis, we have $[v^k]\in N$ which is a right $A_g(V)$-module. Hence $[v^k]\ast_g [a^k]\in N$. Moreover, since $\deg (a^k_{j-1}v^k)=\wt a^k-j+\deg v^k<n$ for any $j\geq 1$, we have $[a^k_{j-1}v^k]\in N$ by the induction hypothesis. Thus $[a^k_{-1}v^k]\in N$ in view of \eqref{eqmiddlestep1}. 
If $r>0$, by \eqref{def:OgM} we have the following equation in $A_g(M)$: 
$$[a^k\circ_g v^k]=[a^k_{-1}v^k]+\sum_{j\geq 1}\binom{\wt a^k-1+\frac{r}{T}}{j} [a^k_{j-1}v^k]=0.$$
 Since $\wt (a^k_{j-1}v^k)<n$ for any $j\geq 1$, we have $[a^k_{-1}v^k]=-\sum_{j\geq 1} \binom{\wt a^k-1+\frac{r}{T}}{j} [a^k_{j-1}v^k]\in N$ by the induction hypothesis. This proves $[M(n)]\ssq N$ and finishes the induction step. Using a similar argument, we can show $	A_g(M)=\sum_{i\in I} A_g(V)\ast_g [u^i]$. Assume $V$ is $C_1$-cofinite and $M$ is $C_1$-cofinite. By Theorem~\ref{thm:AgVN}, $A_g(V)$ is a left (resp. right) Noetherian algebra. By \eqref{AgMgenerate}, $A_g(M)$ is a finitely generated left (resp. right) $A_g(V)$-module. Thus, $A_g(M)$ is left (resp. right) Noetherian as a left (resp. right) $A_g(V)$-module. 

The proof of (2) is similar to the proof of (1) and the proof of Proposition 3.16 in \cite{L99}, we briefly sketch it. Again, we may denote $\sum_{j\in J} A_g(V)\ast_g [w^j]\ast_g A_g(V)$ by $N'$, and use induction on the degree $n$ to show that $[M(n)]\ssq N'$. Assume the conclusion holds for smaller $n$, for $y\in M(n)=W\cap M(n)+\widetilde{C}^g_1(M)\cap M(n)$, we may express it as
$$y=w+\sum_{k=1}^s a^k_{-1}v^k+\sum_{l=1}^t b^l_{0}u^l,\quad w\in W,\ a^k\in V_+\cap V^r,\ b^l\in \oplus_{p\geq 2}V_p\cap V^0,\ v^k,u^l\in M,$$
with $\deg(a^k_{-1}v^k)=\deg (b^l_{0}u^l)=n$ for all $k,l$. By adopting a similar argument as above, we can show $[a^k_{-1}v^k]\in N'$ for all $k$. Moreover, using the facts that 
$$
b\ast_g u-u\ast_g b\equiv \Res_z Y_M(b,z)u (1+z)^{\wt b-1} \pmod{O_g(M)},
$$
for $b\in V^0$ and $u\in M$, and $\deg (b^l_ju^l)<n$ for any $j\geq 1$, we have
\begin{equation}\label{3.6}
[b^l_0u^l]=-\sum_{j\geq 1} \binom{\wt b^l-1}{j} [b^l_j u^l]+[b^l]\ast_g [u^l]-[u^l]\ast_g [b^l]\in N' 
\end{equation}
for all $l$ by the induction hypothesis. Thus $[y]=[w]+\sum_{k=1}^s [a^k_{-1}v^k]+\sum_{l=1}^t [b^l_{0}u^l]\in N'$. 
	\end{proof}
 Using Theorem~\ref{thm:AgMfg}, we can generalize Corollary 3.17 in \cite{L99} and Theorem 3.1 in \cite{H05} about the finiteness of fusion rules under $C_1$-cofinite condition to the $g$-twisted case:

\begin{coro}\label{corofinitenessoffusion}
Let $M^1$ be an untwisted ordinary $V$-module, and $M^2,M^3$ be $g$-twisted ordinary $V$-modules. If the $M^1$ is weakly $C^g_1$-cofinite, then the fusion rule $N\fusion{M^1}{M^2}{M^3}$ is finite. 
	\end{coro}
\begin{proof}
	Since $M^1$ is $C_1$-cofinite implies $M^1$ is weakly $C_1$-cofinite, it suffices to prove the finiteness of $N\fusion{M^1}{M^2}{M^3}$ when $M^1$ is weakly $C_1$-cofinite. 
	The following estimate for the fusion rule was proved by Gao, the author, and Zhu, see \cite{GLZ23} Theorem 6.5: 
	\begin{equation}\label{est:fusion}
	N\fusion{M^1}{M^2}{M^3}\leq \dim (M^3(0)^\ast\otimes_{A_g(V)}A_g(M)\otimes_{A_g(V)}M^2(0))^\ast .\end{equation}
Let $M=W+\widetilde{C}_1(M)$, where $W=\spn\{w^1,\ds, w^n\}$. Then by \eqref{AgMbigenerate}, we have 
\begin{align*}
M^3(0)^\ast\otimes_{A_g(V)}A_g(M)\otimes_{A_g(V)}M^2(0)&=\sum_{j=1}^n M^3(0)^\ast\otimes_{A_g(V)}A_g(V)\ast_g [w^j]\ast_g A_g(V)\otimes_{A_g(V)}M^2(0)\\
&=\sum_{j=1}^n M^3(0)^\ast \otimes_\C \C [w^j]\otimes_\C M^2(0),
\end{align*}
which is finite-dimensional since $M^2(0)$ and $M^3(0)$ are both finite-dimensional by Definition~\ref{def:gtwisted}. Hence $N\fusion{M^1}{M^2}{M^3}$ is finite, in view of \eqref{est:fusion}. 
	\end{proof}

\section{Example of a non-$C_1$-cofinite CFT-type VOA whose Zhu's algebra is non-Noetherian}\label{sec4} 
For the classical non-$C_2$-cofinite CFT-type VOAs, like the vacuum module VOA $V_{\hat{\g}}(\ell,0)$, the Heisenberg VOA $M_{\hat{\h}}(\ell,0)$, and the universal Virasoro VOA $\bar{V}(c,0)$ \cite{FZ92,LL04}, it is well-known that they are $C_1$-cofinite \cite{DLM02}. Although one can construct a non-$C_1$-cofinite VOA by taking infinite direct sum of a $C_1$-cofinite CFT-type VOA, such examples are not of CFT-type. %Natural examples of CFT-type non-$C_1$-cofinite VOAs are quite rare. 

In this Section, we give a natural example of non-$C_1$-cofinite CFT-type VOA based on the idea of the Borel-type subVOAs of a lattice VOA in \cite{Liu24}.  We will also determine its Zhu's algebra and show that it is not Noetherian. 

\subsection{The non-$C_1$-cofinite VOA $V_M$} We refer to \cite{FLM88,DL93} for the general construction of lattice VOAs. Let $L$ be a positive definite even lattice. It was observed by the author in \cite{Liu24} that for any additive submonoid $M\leq L$ with identity element $0$, the subspace 	$V_M:=\bigoplus_{\ga\in M} M_{\hat{\h}}(1,\ga)$ is a CFT-type subVOA of $V_L$, where $\h=\C\otimes_{\Z} L$. See also \cite{DL93}. This follows from the fact that each $M_{\hat{\h}}(1,\ga)$ is a simple current module over the Heisenberg subVOA $M_{\hat{\h}}(1,0)\leq V_L$.

Now let $L=A_2=\Z\al\op \Z\b$ be the type $A_2$ root lattice, where $(\al|\al)=2, (\b|\b)=2,$ and $ (\al|\b)=-1$. Choose a two-cocycle $\epsilon: A_2\times A_2\ra \{\pm 1\}$, where $\epsilon(\al,\al)=1,\epsilon(\b,\b)=1,\epsilon(\al,\b)=1, $ and $\epsilon(\b,\al)=-1$. 
Let $M\leq A_2$ be the following submonoid: 
\begin{equation}\label{def:M}
	M:=\{m\al+n\b: m\geq n\geq 1\}\cup \{0\}.
\end{equation}
$M$ is represented in Figure~\ref{fig1} by the red dots. 

\begin{figure}
	\centering
	\begin{tikzpicture}
	\coordinate (Origin)   at (0,0);
	\coordinate (XAxisMin) at (-8,0);
	\coordinate (XAxisMax) at (8,0);
	\coordinate (YAxisMin) at (0,8);
	\coordinate (YAxisMax) at (0,8);
	%\draw [thin, gray,-latex] (XAxisMin) -- (XAxisMax) node[right] {$x$};% Draw x axis
	%\draw [thin, gray,-latex] (YAxisMin) -- (YAxisMax) node[above] {$z$} ;% Draw y axis
	
	\clip (-4.6,-2.3) rectangle (5.5,6); % Clips the picture...
	\begin{scope} %<- added
	\pgftransformcm{1}{0}{1/2}{sqrt(3)/2}{\pgfpoint{0cm}{0cm}} 
	\coordinate (Bone) at (0,2);
	\coordinate (Btwo) at (2,-2);
	\draw[style=help lines,dashed] (-7,-6) grid[step=2cm] (6,6);
	\foreach \x in {-4,-3,...,4}{% Two indices running over each
		\foreach \y in {1,...,4}{% node on the grid we have drawn 
			\coordinate (Dot\x\y) at (2*\x,2*\y);
			\node[draw,circle,inner sep=2pt,fill] at (Dot\x\y) {};
		}
	}
	\foreach \x in {1,...,4}{% Two indices running over each
		\foreach \y in {0}{% node on the grid we have drawn 
			\coordinate (Dot\x\y) at (2*\x,2*\y);
			\node[draw,circle,inner sep=2pt,fill] at (Dot\x\y) {};
		}
	}
	\foreach \x in {0,...,4}{% Two indices running over each
		\foreach \y in {1,...,4}{% node on the grid we have drawn 
			\coordinate (Dot\x\y) at (2*\x,2*\y);
			\node[draw,red,circle,inner sep=2pt,fill] at (Dot\x\y) {};
		}
	}
	\foreach \x in {-1,...,-4}{% Two indices running over each
		\foreach \y in {0}{% node on the grid we have drawn 
			\coordinate (Dot\x\y) at (2*\x,2*\y);
			\node[draw,circle,inner sep=2pt,fill] at (Dot\x\y) {};
		}
	}
	
	\foreach \x in {-4,-3,...,4}{% Two indices running over each
		\foreach \y in {-1,...,-4}{% node on the grid we have drawn 
			\coordinate (Dot\x\y) at (2*\x,2*\y);
			\node[draw,circle,inner sep=2pt,fill] at (Dot\x\y) {};
		}
	}
	\foreach \x in {0}{% Two indices running over each
		\foreach \y in {0}{% node on the grid we have drawn 
			\coordinate (Dot\x\y) at (2*\x,2*\y);
			\node[draw,red,circle,inner sep=2pt,fill] at (Dot\x\y) {};
			\node [xshift=0cm, yshift=-0.4cm]{$0$};
		}
	}
	
	\draw [thick,-latex,red] (Origin) 
	-- (Bone) node [above right]  {$\al+\b$};
	\draw [thick,-latex] (Origin)
	-- (Btwo) node [below right]  {$-\b$};
	\draw [thick,-latex] (Origin) node [xshift=0cm, yshift=-0.3cm]{}
	-- ($(Bone)+(Btwo)$) node [below right] {$\al$};
	\draw [thick,-latex] (Origin)
	-- ($-1*(Bone)-1*(Btwo)$) node [below left] {$-\al$};
	\draw [thick,-latex] (Origin)
	-- ($-1*(Btwo)$) coordinate (B3) node [above left] {$\b$};
	\draw [thick,-latex] (Origin)
	-- ($-1*(Bone)$) node [below left]  {$-\al-\b$};
	% all the following is added
	\end{scope} 
	\begin{scope}
	\pgftransformcm{1}{0}{-1/2}{sqrt(3)/2}{\pgfpoint{0cm}{0cm}} 
	\draw[style=help lines,dashed] (-6,-6) grid[step=2cm] (6,6);
	\end{scope}
	\end{tikzpicture}
	\caption{ \label{fig1}}
\end{figure}

Consider the subVOA of the lattice VOA $V_{A_2}$ associated to $M$: 
\begin{equation}\label{structureV_M}
V_M=M_{\hat{\h}}(1,0)\op\bigoplus_{m\geq n\geq 1} M_{\hat{\h}}(1,m\al+n\b). 
\end{equation} 
 Then $V_M$ is of CFT-type. In the rest of this Section, we fix the VOA $V_M$ as in \eqref{structureV_M}. We will show that $V_M$ is not $C_1$-cofinite, and $A(V_M)$ is not Noetherian.

\begin{lm}\label{lmforC_1M}
For any $m\geq 1$, we have $e^{m\al+\b}\notin C_1(V_M)$. 
	\end{lm}
\begin{proof}
In view of \eqref{structureV_M} and \eqref{def:C_1}, we can express $C_1(V_M)$ as follows: 
\begin{align*}
C_1(V_M)&=\spn\{u_{-1}v: u\in M_{\hat{\h}}(1,\ga)\cap (V_M)_+, v\in M_{\hat{\h}}(1,\ga')\cap(V_M)_+, \ga,\ga'\in M \}\\
&\ +\spn\{L(-1)w:w\in M_{\hat{\h}}(1,\theta)\cap(V_M)_+, \theta\in M \}.
\end{align*}
Suppose $e^{m\al+\b}\in C_1(V_M)$ for some $m\geq 1$. 
Note that $u_{-1}v\in  M_{\hat{\h}}(1,\ga+\ga')$ for $u\in  M_{\hat{\h}}(1,\ga)$ and $v\in  M_{\hat{\h}}(1,\ga')$, and $L(-1) M_{\hat{\h}}(1,\theta)\ssq  M_{\hat{\h}}(1,\theta)$. Moreover, if $\ga,\ga'$ are nonzero elements in $M$ \eqref{def:M}, then $\ga+\ga'\neq m\al+\b$. Since the Heisenberg modules $M_{\hat{\h}}(1,m\al+n\b)$ are in direct sum in \eqref{structureV_M}, and $C_1(V_M)$ is a graded subspace of $V_M$, it follows that $e^{m\al+\b}$ must be contained in $W_1+W_2+W_3\subset C_1(V_M)$, where
\begin{align*}
W_1&= \spn\{u_{-1}v: u\in  M_{\hat{\h}}(1,m\al+\b)\cap (V_M)_+, v\in M_{\hat{\h}}(1,0)\cap (V_M)_+  \},\\
W_2&=\spn\{u'_{-1}v': u'\in M_{\hat{\h}}(1,0)\cap (V_M)_+,\ v'\in M_{\hat{\h}}(1,m\al+\b)\cap (V_M)_+ \},\\
W_3&=\spn\{L(-1)w:  w\in M_{\hat{\h}}(1,m\al+\b)\cap (V_M)_+\}.
\end{align*}
Note that $W_1,W_2\subset  M_{\hat{\h}}(1,m\al+\b)=\bigoplus_{k=0}^\infty M_{\hat{\h}}(1,m\al+\b)_{(m^2-m+1)+k}$, where $m^2-m+1=\wt (e^{m\al+\b})$. For homogeneous elements $u\in  M_{\hat{\h}}(1,m\al+\b)\cap (V_M)_+$ and $ v\in M_{\hat{\h}}(1,0)\cap (V_M)_+ $, since $\wt v>0$, we must have 
$u_{-1}v\in \sum_{k=1}^\infty M_{\hat{\h}}(1,m\al+\b)_{(m^2-m+1)+k}$ as $\wt (u_{-1}v)=\wt u+\wt v>\wt u\geq m^2-m+1$. This shows $W_1,W_2\ssq  \sum_{k=1}^\infty M_{\hat{\h}}(1,m\al+\b)_{(m^2-m+1)+k}$. On the other hand, since $\wt (L(-1)w)=\wt w+1$, it is clear that $W_3\ssq \sum_{k=1}^\infty M_{\hat{\h}}(1,m\al+\b)_{(m^2-m+1)+k}$. Then we have 
$$e^{m\al+\b}\in M_{\hat{\h}}(1,m\al+\b)_{m^2-m+1}\cap  \sum_{k=1}^\infty M_{\hat{\h}}(1,m\al+\b)_{(m^2-m+1)+k}=0,$$
which is a contradiction. Thus, $e^{m\al+\b}\notin C_1(V_M)$ for any $m\geq 1$. 
	\end{proof}

%By abuse of notations, we use the symbol $[u]$ (similar to $[a]$ in $A_g(V)$) to denote the equivalent class of $u \in V_M$ in the quotient space $V_M/C_1(V_M)$. 

\begin{thm}\label{thm:V_MisnotC1}
	 $V_M/C_1(V_M)$ has a basis $\{\vac+C_1(V_M),e^{m\al+\b}+C_1(V_M):m\geq 1\}$. In particular, the CFT-type VOA $V_M$ is {\em  not} $C_1$-cofinite. 
	\end{thm}
	\begin{proof}
	Since $(m\al+n\b|m\al+n\b)/2=m^2-mn+n^2\geq 1$ for all $m\geq n\geq 1$, we have $M_{\hat{\h}}(1,m\al+n\b)\ssq (V_M)_+$ for any such a pair of $m,n$. Also note that $a_{-n}b\in C_1(V_M)$ for any $a,b\in (V_M)_+$ and $n\geq 1$ since $C_1(V)\supset C_2(V)\supset C_3(V)\supset \ds $, see \cite{L99}.

First, we show that $e^{m\al+n\b}\in C_1(V_M)$ for any $m\geq n\geq 2$. Indeed, for any $m\geq n\geq 1$, since $(\al+\b|m\al+n\b)=m+n\geq 2$, by the definition of lattice vertex operators, we have
		\begin{align*}
		e^{\al+\b}_{-m-n-1}e^{m\al+n\b}&=\Res_z E^-(-\al-\b,z)E^+(-\al-\b,z)e_{\al+\b}z^{\al+\b}e^{m\al+n\b}\\
		&=(-1)^m e^{(m+1)\al+(n+1)\b}\in C_1(V_M).
		\end{align*}
		Hence $e^{(m+1)\al+(n+1)\b}\in C_1(V_M)$ for any $(m+1)\geq (n+1)\geq 2$. This proves $e^{m\al+n\b}\in C_1(V_M)$ for any $m\geq n\geq 2$. Since $h(-n)C_1(V_M)\ssq C_1(V_M)$ for any $h\in \h$ and $n\geq 1$, we have $M_{\hat{\h}}(1,m\al+n\b)\ssq C_1(V_M)$ for any $m\geq n\geq 2$. Then by the decomposition \eqref{structureV_M}, we have
	\begin{align*}
	V_M/C_1(V_M)&=\left(M_{\hat{\h}}(1,0)+\sum_{m\geq 1} M_{\hat{\h}}(1,m\al+\b)\right)+C_1(V_M)\\
	&= \spn\{\vac+C_1(V_M), e^{m\al+\b}+C_1(V_M): m\geq 1 \}. 
	\end{align*}
It remains to show $\{\vac+C_1(V_M), e^{m\al+\b}+C_1(V_M): m\geq 1 \}$ are linearly independent. 

Note that $\vac+C_1(V_M)\neq 0$ in view of \eqref{def:C_1}. By Lemma~\ref{lmforC_1M}, $e^{m\al+\b}+C_1(V_M)\neq 0$ for any $m\geq 1.$ Since $L(0)(u_{-1}v)=u_{-1}L(0)v+(L(0)u)_{-1}v$ and $L(0)L(-1)w=L(-1)L(0)w$, we have $L(0)C_1(V_M)\ssq C_1(V_M)$. Hence $L(0): V_M/C_1(V_M)\ra V_M/C_1(V_M), a+C_1(V_M)\mapsto L(0)a+C_1(V_M)$ is a well-defined linear map. Since $L(0)e^{m\al+\b}=(m^2-m+1) e^{m\al+\b}$, it follows that $\vac+C_1(V_M), e^{\al+\b}+C_1(V_M),e^{2\al+\b}+C_1(V_M),e^{3\al+\b}+C_1(V_M),\dots $  are eigenvectors of $L(0)$ of distinct eigenvalues. Thus, $\{\vac+C_1(V_M),e^{m\al+\b}+C_1(V_M):m\geq 1\}$ is a basis of  $V_M/C_1(V_M)$. 
		\end{proof}
	
\begin{remark}
From the proofs of Lemma~\ref{lmforC_1M} and Theorem~\ref{thm:V_MisnotC1}, we see that 
the essential reason why $V_M$ is not $C_1$-cofinite is that the chain of lattice points $\{\al+\b, 2\al+\b,3\al+\b,\ds \}$, which is the first horizontal row of red dots in Figure~\ref{fig1}, cannot be generated by finitely many points in the submonoid $M$. Using this idea, one can construct many examples of non-$C_1$-cofinite CFT-type VOAs inside a lattice VOA.
\end{remark}

\subsection{Non-Noetherianity of the Zhu's algebra of $V_M$} We determine the (untwisted) Zhu's algebra of $V_M$ \eqref{structureV_M} based on a similar method as in \cite{Liu24}. In particular, we will see that $A(V_M)$ is {\em not} Noetherian. Hence our example $V_M$ in this Section verifies Theorem~\ref{thm:AgVN} when $g=\Id_V$. 

\subsubsection{A spanning set of $O(V_M)$}

\begin{lm}\label{lm4.3}
For any $m\geq n\geq 2$, we have $e^{m\al+n\b}\in O(V_M)$. 
\end{lm}
\begin{proof}
Similar to the proof of Theorem~\ref{thm:V_MisnotC1}, for any $m\geq n\geq 1$, since $(\al+\b|m\al+n\b)=m+n\geq 2$, we have the following formula for any $k\geq 0$: 
\begin{align*}
e^{\al+\b}_{-k-1}e^{m\al+n\b}&=\Res_z (-1)^mE^-(-\al-\b,z) z^{(\al+\b|m\al+n\b)-k-1} e^{(m+1)\al+(n+1)\b}\\
&=\begin{cases}(-1)^m e^{(m+1)\al+(n+1)\b}&\mathrm{if}\ k=m+n,\\
0&\mathrm{if}\ k<m+n. \end{cases} 
\end{align*}
 Then by \eqref{propertyofOgM} with $M=V$ and $g=\Id_V$, noting that $\wt (e^{\al+\b})=2$,  we have 
\begin{align*}
\Res_z Y(e^{\al+\b},z) e^{m\al+n\b}\frac{(1+z)^2}{z^{m+n+1}}&=e^{\al+\b}_{-m-n-1}e^{m\al+n\b}+2e^{\al+\b}_{-m-n}e^{m\al+n\b}+e^{\al+\b}_{-m-n+1}e^{m\al+n\b}\\
&=(-1)^m e^{(m+1)\al+(n+1)\b}\equiv 0\pmod{O(V_M)}. 
\end{align*}
Thus, $e^{(m+1)\al+(n+1)\b}\in O(V_M)$ for any $(m+1)\geq (n+1)\geq 2$. 
	\end{proof}

Let $O$ be the subspace of $V_M$ spanned by the following elements:
\begin{equation}\label{def:O}
\begin{cases}
&h(-n-2)u+h(-n-1)u, \qquad u\in V_M,\ \mathrm{and}\ n\geq 0,\\
&m\al(-1)v+\b(-1)v+(m^2-m+1)v,\qquad v\in M_{\hat{\h}}(1,m\al+\b),\ m\geq 1,\\
&M_{\hat{\h}}(1,m\al+n\b),\qquad m\geq n\geq 2.
\end{cases}\end{equation}

\begin{lm}\label{prop:Ois asubsetofOVM}
We have $O\ssq O(V_M)$ as subspace of $V_M$. 
	\end{lm}
\begin{proof}
Clearly, $h(-n-2)u+h(-n-1)u\in O(V_M)$ for any $u\in V_M$. Let $m\geq 1$. Note that $e^{m\al+\b}_{-2}\vac= \Res_z E^-(-\al-\b,z) z^{-2}e^{m\al+\b}=m\al(-1)e^{m\al+\b}+\b(-1)e^{m\al+\b}$. Since $\wt (e^{m\al+\b})=m^2-m+1$, we have 
\begin{align*}
e^{m\al+\b}\circ \vac&=\Res_z Y(e^{m\al+\b},z)\vac\frac{(1+z)^{m^2-m+1}}{z^2}= e^{m\al+\b}_{-2}\vac+ (m^2-m+1)e^{m\al+\b}_{-1}\vac\\
&=m\al(-1)e^{m\al+\b}+\b(-1)e^{m\al+\b}+(m^2-m+1)e^{m\al+\b}\equiv 0\pmod{O(V_M)}.
\end{align*}
Since $\al(-1)$ and $\b(-1)$ commute with $h(-n)$, for any $h\in \h$ and $n\geq 1$, we have $m\al(-1)v+\b(-1)v+(m^2-m+1)v\in O(V_M)$, for any $v=h^1(-n_1)\ds h^r(-n_r)e^{m\al+\b}\in M_{\hat{\h}}(1,0)$.  

Finally, for $m\geq m\geq 2$, let $w=h^1(-n_1-1)\ds h^r(-n_r-1)e^{m\al+\b}$ be a spanning element of $M_{\hat{\h}}(1,m\al+n\b)$, where $n_1\geq \ds \geq n_r\geq 0$. Since $h(-n-1)v\equiv (-1)^n v\ast (h(-1)\vac)\pmod{O(V_M)}$ for any $h\in \h$ and $v\in V_M$ by \eqref{rightmodact}, we have 
$$
w\equiv (-1)^{n_1+\ds +n_r} e^{m\al+\b}\ast (h^1(-1)\vac)\ast \ds \ast (h^r(-1)\vac)\pmod{O(V_M)}.
$$
Moreover, by Theorem 2.1.1 in \cite{Z96}, $O(V_M)\ast V_M\ssq O(V_M)$. It follows from Lemma~\ref{lm4.3} that $w\in O(V_M)$. Hence $M_{\hat{\h}}(1,m\al+n\b)\ssq O(V_M)$ for any $m\geq m\geq 2$.
	\end{proof}

Conversely, we want to show $O(V_M)\ssq O$. By \eqref{structureV_M}, it suffices to show 
$M_{\hat{\h}}(1,\ga)\circ M_{\hat{\h}}(1,\ga')\ssq O,$ for any $\ga,\ga'\in M$. If $\ga=m\al+n\b$ and $\ga'=m'\al+n'\b$, where $m\geq n\geq 1$ and $m'\geq n'\geq 1$, then $\ga+\ga'=(m+m')\al+(n+n')\b$, with $m+m'\geq n+n'\geq 2$. By \eqref{def:O} we have 
$$M_{\hat{\h}}(1,\ga)\circ M_{\hat{\h}}(1,\ga')\ssq M_{\hat{\h}}(1,\ga+\ga')\subset O.$$
Moreover, if $m\geq n\geq 2$, we also have $M_{\hat{\h}}(1,0)\circ M_{\hat{\h}}(1,m\al+n\b)\subset O$ and $M_{\hat{\h}}(1,m\al+n\b)\circ M_{\hat{\h}}(1,0)\subset O$. Hence we only need to show 
\begin{equation}\label{4.4}
M_{\hat{\h}}(1,0)\circ M_{\hat{\h}}(1,m\al+\b)\subset O,\quad \mathrm{and}\quad 
M_{\hat{\h}}(1,m\al+\b)\circ M_{\hat{\h}}(1,0)\subset O,
\end{equation}
for any $m\geq 1$. The proof of \eqref{4.4} is a slight modification of the induction process in Section 3.2 in \cite{Liu24}, we omit the details. In conclusion, we have the following: 

\begin{prop}\label{propO}
Let $O$ be the subspace of $V_M$ spanned by elements in \eqref{def:O}. Then $O=O(V_M)$. 
	\end{prop}

\subsubsection{Structure of $A(V_M)$} Consider the associative algebra 
\begin{equation}
A_M=\C\<x,y,z_1,z_2,\ds\>/R,
\end{equation}
where $\C\<x,y,z_1,z_2,\ds\>$ is the tensor algebra on infinitely many generators $x,y,z_1,z_2,\ds $, and $R$ is the two-sided ideal generated by the following elements: 
\begin{equation}\label{relation}
\begin{aligned}
&xy-yx,\quad  z_m(mx+y)+(m^2-m+1)z_m,\ m\geq 1,\quad  z_iz_j,\ i,j\geq 1,\\
& xz_m-z_mx-(2m-1)z_m,\quad yz_m-z_my-(2-m)z_m,\ m\geq 1.
\end{aligned}
\end{equation}

It is clear that $A_M$ has the following subspace decomposition:
\begin{equation}\label{decofA_M}
A_M=\C[x,y]\op \left(\bigoplus_{m=1}^\infty z_m\C[y]\right),
\end{equation}
where $z_m\C[y]$ is a vector space with basis $\{z_m,z_my,z_my^2,\ds\}$, and we use the same symbols $x,y,z_m$ to denote their equivalent classes in the quotient space. 

\begin{thm}\label{thm4.6}
	Define an algebra homomorphism $F: \C\<x,y,z_1,z_2,\ds\>\ra A(V_M)$ by letting 
	\begin{equation}\label{4.7}
	F(x):=[\al(-1)\vac],\ F(y):=[\b(-1)\vac],\ F(z_m)=[e^{m\al+\b}],\ m\geq 1.
	\end{equation}
	Then $F$ factors through $A_M$ and induces an isomorphism $F:A_M\ra A(V_M)$. 
	\end{thm}
\begin{proof}
%The fact that $F$ factors through $A_M$ i.e., $F(R)=0$ is an easy consequence of \eqref{def:O}, see also Proposition 3.1 in \cite{Liu24}. We omit the details. 
We first show that $F(R)=0$. Indeed, by \eqref{commprodAgV}, \eqref{def:O}, and Lemma~\ref{prop:Ois asubsetofOVM}, we have 
\begin{align*}
&F(xy-yx)=[\al(-1)\vac]\ast[\b(-1)\vac]-[\b(-1)\vac]\ast[\al(-1)\vac]\\
& =[\al(0)\b(-1)\vac]=0;\\
&F(z_m(mx+y)+(m^2-m+1)z_m)\\
&= [e^{m\al+\b}]\ast [m\al(-1)\vac]+[e^{m\al+\b}]\ast [\b(-1)\vac]+(m^2-m+1)[e^{m\al+\b}] \\
&=[m\al(-1)e^{m\al+\b}+\b(-1)e^{m\al+\b}+(m^2-m+1)e^{m\al+\b}]=0;\\
&F(xz_m-z_mx-(2m-1)z_m)=[\al(-1)\vac]\ast [e^{m\al+\b}]-[e^{m\al+\b}]\ast [\al(-1)\vac]-(2m-1)[e^{m\al+\b}]\\
&=[\al(0)e^{m\al+\b}-(2m-1)e^{m\al+\b}]=0;\\
& F(yz_m-z_my-(2-m)z_m)= [\b(-1)\vac]\ast [e^{m\al+\b}]-[e^{m\al+\b}]\ast [\b(-1)\vac]-(2-m)[e^{m\al+\b}]\\
&=[\b(0)e^{m\al+\b}-(2-m)e^{m\al+\b}]=0.
\end{align*}
Moreover, since $e^{i\al+\b}\ast e^{j\al+\b}\in M_{\hat{\h}}(1,(i+j)\al+2\b)\subset O(V_M)$ in view of  Proposition~\ref{prop:Ois asubsetofOVM}, it follows that 
$F(z_iz_j)=[e^{i\al+\b}\ast e^{i\al+\b}]=0,$
for any $i,j\geq 1$. This shows $F$ factors though $A_M$. To show $F$ is an isomorphism, we construct an inverse map of $F$. Similar to the proof of Theorem 4.11 in \cite{Liu24}, we first define a linear map 
\begin{equation}\label{linearmap}
\bar{(\cdot)}: \h=\C\al \op \C \b\ra A_M, \quad h=\la\al+\mu \b\mapsto \bar{h}=\la x+\mu y,\quad \la,\mu\in\C.
\end{equation}
Then we define a linear map $G: V_M=M_{\hat{\h}}(1,0)\op\bigoplus_{m\geq n\geq 1} M_{\hat{\h}}(1,m\al+n\b) \ra A_M$ as follows: 
\begin{align}
h^1(-n_1-1)\ds h^r(-n_r-1)\vac&\mapsto (-1)^{n_1+\ds +n_r} \overline{h^r}\cdot  \overline{h^{r-1}}\cdot \ds \cdot \overline{h^1}, \label{4.9}\\
h^1(-n_1-1)\ds h^r(-n_r-1)e^{m\al+\b}&\mapsto (-1)^{n_1+\ds +n_r} z_m\cdot \overline{h^r}\cdot  \overline{h^{r-1}}\cdot \ds \cdot \overline{h^1},\quad m\geq 1, \label{4.10}\\
M_{\hat{\h}}(1,m\al+n\b)&\mapsto 0,\quad m\geq n\geq 2,\label{4.11}
\end{align}
where $n_1\geq \ds\geq n_r\geq 0$, and $\overline{h^i}$ is the image of $h^i$ in $A_M$ under \eqref{linearmap}. We claim that $G(O)=0$. 

Indeed, for any $u\in V_M$, $h\in \h$, and $n\geq 0$, it is clear from \eqref{4.9}--\eqref{4.11} that $$G(h(-n-2)u+h(-n-1)u)=(-1)^{n+1}G(u)\cdot \overline{h}+(-1)^n G(u)\cdot \overline{h}=0.$$
Let $v=h^1(-n_1-1)\ds h^r(-n_r-1)e^{m\al+\b}$ be a spanning element of $ M_{\hat{\h}}(1,m\al+\b)$ with $m\geq 1$, then by \eqref{linearmap}, \eqref{4.10}, and \eqref{relation}, we have 
\begin{align*}
&G(m\al(-1)v+\b(-1)v+(m^2-m+1)v)\\
&=m  (-1)^{n_1+\ds +n_r} z_m\cdot \overline{h^r}\cdot  \ds \cdot \overline{h^1}\cdot \overline{\al} +  (-1)^{n_1+\ds +n_r} z_m\cdot \overline{h^r}\cdot \ds \cdot \overline{h^1}\cdot \overline{\b}+ (m^2-m+1)z_m\cdot \overline{h^r}\cdot  \ds \cdot \overline{h^1}\\
&=(-1)^{n_1+\ds+n_r}\left( z_m(mx+y)+(m^2-m+1)z_m \right)\cdot  \overline{h^r}\cdot  \ds \cdot \overline{h^1}\\
&=0.
\end{align*}
Finally, by \eqref{4.11} we have $G(M_{\hat{\h}}(1,m\al+n\b))=0$ for any $m\geq n\geq 2$. Thus, we have $G(O(V_M))=0$ by \eqref{def:O} and Proposition~\ref{propO}, and $G$ induces a well-defined linear map $G:A(V_M)=V_M/O(V_M)\ra A_M$, such that 
\begin{equation}\label{4.12}
G([\al(-1)\vac])=x,\quad G([\b(-1)\vac])=y,\quad G([e^{m\al+\b}])=z_m,\ m\geq 1,
\end{equation}
in view of \eqref{4.9}--\eqref{4.11}. By \eqref{4.7} and \eqref{4.12}, it is clear that $G$ is an inverse of $F:A_M\ra A(V_M)$. Hence $A_M\cong A(V_M)$ as associative algebras. 
	\end{proof}

\begin{coro}\label{corononnoehterian}
The untwisted Zhu's algebra $A(V_M)$ is {\em not} Noetherian. 
	\end{coro}
\begin{proof}
	By Theorem~\ref{thm4.6} and \eqref{decofA_M}, we have an isomorphism 
	$$A(V_M)\cong A_M= \C[x,y]\op \left(\bigoplus_{m=1}^\infty z_m\C[y]\right)=\C[x,y]\op J,$$
	where $J=\bigoplus_{m=1}^\infty z_m\C[y]$. By \eqref{relation}, it is clear that $J$ is a two-sided ideal of $A_M$. Suppose $J$ can be generated by finitely many elements $w_1,\ds,w_k\in J$. There must exist an index $N>0$ s.t. 
	$w_1,\ds ,w_k\in \bigoplus_{m=1}^N z_m \C[y]$. But it follows from \eqref{relation} that 
	$$A_M\cdot \left(\bigoplus_{m=1}^N z_m \C[y]\right)\cdot A_M\ssq \bigoplus_{m=1}^N z_m \C[y],$$
	since $z_iz_j=0$ for all $i,j\geq 1$. Then we have $J\ssq \bigoplus_{m=1}^N z_m \C[y]$, which is a contradiction. Therefore, $A(V_M)$ has a two-sided ideal $J$ that is not finitely generated. This shows $A(V_M)$ is neither left nor right Noetherian.  
	\end{proof}

\begin{remark}
	Using Corollary~\ref{corononnoehterian} and Theorem~\ref{thm:AgVN} with $g=\Id_V$, we have an alternative proof of the fact that $V_M$ is not $C_1$-cofinite without finding a basis of $V_M/C_1(V_M)$ as in Theorem~\ref{thm:V_MisnotC1}. 
	\end{remark}

%-------------------------------------------------------------------------------------------------------

\section{Finiteness of $g$-twisted higher Zhu's algebra}\label{Sec4}

In this Section, using a higher order analog of the epimorphism \eqref{varphi}, we prove that the $g$-twisted higher Zhu's algebra $A_{g,n}(V)$ constructed by Dong, Li, and Mason in \cite{DLM98(2)} and its bimodule $A_{g,n}(M)$ constructed by Jiang and Jiao in \cite{JJ16} are finite-dimensional if $V$ is $C_2$-cofinite, which generalizes Miyamoto's result on finiteness of $A_{n}(V)$ and Buhl's result on finiteness of $A_n(M)$ under the $C_2$-cofinite condition in \cite{M04} to the $g$-twisted case. 

\subsection{Shifted Level filtration on $A_{g,n}(V)$}\label{sec4.1}
 First, we recall the definition of $A_{g,n}(V)$ in \cite{DLM98(2)}. Fix a rational number $n=l+\frac{i}{T}\in \frac{1}{T}\Z$, where $l\in \N$ and $0\leq i\leq T-1$ are uniquely determined by $n$. 
 
 For $a\in V^r$ with $0\leq r\leq T-1$, and $b\in V$, define
 \begin{equation}\label{def:OgnV}
a\circ_{g,n} b:=\Res_z Y(a,z)b\frac{(1+z)^{\wt a-1+\delta_{i}(r)+l+r/T}}{z^{2l+\delta_i(r)+\delta_i(T-r)}},\quad \mathrm{where}\quad 
\delta_i(r)=\begin{cases}
1&\mathrm{if}\ r\leq i\\
0&\mathrm{if}\ r>i
\end{cases},
 \end{equation}
 and set $\delta_i(T)=1$. Let $O_{g,n}(V)$ be the subspace of $V$ spanned by all $a\circ _{g,n}b$ and $L(-1)c+L(0)c$, and let $A_{g,n}(V):=V/O_{g,n}(V)$. Define 
 \begin{equation}\label{AgnVprod}
 a\ast_{g,n}b:=\begin{cases}\sum_{m=0}^l (-1)^m\binom{m+l}{l} \Res_z Y(a,z)b (1+z)^{\wt a+l}/z^{l+m+1} &\mathrm{if}\ a\in V^0,\\
 0&\mathrm{if}\ a\in V^r,\ r>0.
 \end{cases}
 \end{equation}
By Theorem 2.4 in \cite{DLM98(2)}, $A_{g,n}(V)$ is an associative algebra with respect to \eqref{AgnVprod}. Again, we denote the equivalent class of an element $a\in V$ in $A_{g,n}(V)$ by $[a]$.

For the rest of this paper, we fix the rational number $n=l+\frac{i}{T}$. The usual level filtration \eqref{levelfilAgV} cannot give us a desirable higher order analog of the epimorphism \eqref{varphi}. So we introduce a new level filtration on $A_{g,n}(V)$ as follows: 
For $p\geq 2l$, let
\begin{equation}\label{levelfilonAgnV}
F_pA_{g,n}(V):=\spn\{[a]: a\in V\ \mathrm{homogeneous},\ \wt a+2l\leq p \}. 
\end{equation} 
For $p<2l$, let $F_{p}A_{g,n}(V):=0$. Clearly we have 
\begin{equation}\label{filAgnV}
F_{2l}A_{g,n}(V)\ssq F_{2l+1}A_{g,n}(V)\ssq \ds,\quad \mathrm{and}\quad A_{g,n}(V)=\bigcup_{p=2l}^\infty F_{p}A_{g,n}(V).
\end{equation}

\begin{lm}\label{lmgrAgnV}
$A_{g,n}(V)$ is a filtered algebra with respect to the filtration \eqref{levelfilonAgnV}. The product on the associated graded algebra $\gr A_{g,n}(V)=\bigoplus_{p=2l}^\infty F_{p}A_{g,n}(V)/F_{p-1}A_{g,n}(V)$ is given by 
\begin{equation}\label{prodgrAgnV}
\left([a]+F_{p-1}A_{g,n}(V)\right)\ast _{g,n} \left([b]+F_{q-1}A_{g,n}(V)\right)=\begin{cases}(-1)^l\binom{2l}{l}[a_{-2l-1}b]+F_{p+q-1}A_{g,n}(V)&\mathrm{if}\ r=0\\
0&\mathrm{if}\ r>0
\end{cases},
\end{equation}
where $a\in V^r,b\in V$ are homogeneous, with $\wt a+2l\leq p$ and $\wt b+2l\leq q$, and $p,q\geq 2l$. 
	\end{lm}
\begin{proof}
Let $a\in V^0$ and $b\in V$ be homogeneous elements such that $\wt a+2l\leq p$ and $\wt b+2l\leq q$. By the definitions of product on $A_{g,n}(V)$ \eqref{AgnVprod}, we have 
\begin{equation}\label{usefuleq}
[a\ast_{g,n}b]=\sum_{m=0}^l \sum_{j\geq 0} (-1)^{m} \binom{m+l}{l} \binom{\wt a+l}{j} [a_{j-l-m-1}b]\in F_{p+q}A_{g,n}(V)
\end{equation}
since $\wt (a_{j-l-m-1}b)+2l=\wt a-j+l+m+\wt b+2l\leq (\wt a+2l)+(\wt b+2l)\leq p+q$, for any $j\geq 0$ and $0\leq m\leq l$. Hence $F_{p}A_{g,n}(V)\ast_{g,n} F_qA_{g,n}(V)\ssq F_{p+q}A_{g,n}(V)$, and so $A_{g,n}(V)$ is a filtered algebra. Moreover, by \eqref{levelfilonAgnV} and the equality about weight above, we have $[a_{j-l-m-1}b]\in F_{p+q-1}A_{g,n}(V)$ unless $j=0$ and $m=l$. It follows from \eqref{usefuleq} that 
$$[a\ast_{g,n}b]+F_{p+q-1}A_{g,n}(V)=(-1)^l\binom{2l}{l} [a_{-2l-1}b]+F_{p+q-1}A_{g,n}(V).$$
This proves \eqref{prodgrAgnV}.
	\end{proof}

\begin{remark}\label{commofgrAgnV}
	By Lemma 2.2 in \cite{DLM98(2)}, for any $a\in V^0$ and $b\in V$, we have 
	\begin{equation}
	a\ast_{g,n}b-b\ast_{g,n}a\equiv \Res_z Y(a,z)b (1+z)^{\wt a-1}=\sum_{j\geq 0} \binom{\wt a-1}{j} a_jb \pmod{O_{g,n}(V)}.
	\end{equation}
	Since $\wt (a_jb)+2l\leq (\wt a+2l)+(\wt b+2l)-1$, it follows that $\gr A_{g,n}(V)$ is a commutative graded algebra with respect to the product \eqref{prodgrAgnV}. However, unlike $\gr A_g(V)$ in the previous sections, if $l\geq 1$, the element $[\vac]+F_{2l-1}A_{g,n}(V)$ is {\bf not} the unit element of $\gr A_{g,n}(V)$. 
	\end{remark}

\subsection{Finiteness of $A_{g,n}(V)$}\label{sec4.2}

 Recall that $n=l+\frac{i}{T}$, where $l\in \N$. Consider the following higher level generalization of Zhu's $C_2$-algebra $R_2(V)$: 
\begin{equation}\label{def:R2l+2V}
R_{2l+2}(V):=V/C_{2l+2}(V),\quad \mathrm{where}\quad C_{2l+2}(V)=\spn\{a_{-2l-2}b: a,b\in V\}.
\end{equation}
\begin{lm}\label{lmR2l+2V}
	$R_{2l+2}(V)$ is a graded associative algebra with respect to the product 
	\begin{equation}\label{prodR2l+2V}
	(a+C_{2l+2}(V))\cdot (b+C_{2l+2}(V))=(-1)^l\binom{2l}{l}a_{-2l-1}b+C_{2l+2}(V),
	\end{equation}
	and the grading 
	\begin{equation}\label{gradingR2l+2V}
	R_{2l+2}(V)=\bigoplus_{p=2l}^\infty R_{2l+2}(V)_p,\quad R_{2l+2}(V)_p=\spn\{a+C_{2l+2}(V):\wt a+2l= p\}.
	\end{equation} 
It is commutative if and only if $\sum_{j=1}^{2l} (-1)^j (b_{j-2l-1}a)_{-1-j}\vac\in C_{2l+2}(V)$, for any $a,b\in V$. 
	\end{lm}
\begin{proof}
	For $a,b,c\in V$, by the Jacobi identity of VOA, we have 
	\begin{align*}
&(a_{-2l-1}b)_{-2l-1}c-a_{-2l-1}(b_{-2l-1}c)\\
&=\sum_{j\geq 1}\binom{-2l-1}{j}(-1)^j a_{-2l-1-j}b_{-2l-1+j}c-\sum_{j\geq 0}\binom{-2l-1}{j} (-1)^{-2l-1+j} b_{-4l-2-j}a_{j}c\\
&\equiv 0\pmod{C_{2l+2}(V)}.
	\end{align*}
Hence the product \eqref{prodR2l+2V} is associative since the coefficient $(-1)^l \binom{2l}{l}$ does not depend on $a$ and $b$. By \eqref{def:R2l+2V}, $C_{2l+2}(V)$ is spanned by homogeneous elements, hence $R_{2l+2}(V)$ is a graded algebra. Since $R_{2l+2}(V)_p=(V_{p-2l}+C_{2l+2}(V))/C_{2l+2}(V)$ for any $p\geq 2l$ and $V=\bigoplus_{p=2l}^\infty V_{p-2l}$, we have $	R_{2l+2}(V)=\bigoplus_{p=2l}^\infty R_{2l+2}(V)_p$. It is clear that $R_{2l+2}(V)_{p}\cdot R_{2l+2}(V)_q\ssq R_{2l+2}(V)_{p+q}$ for any $p,q\geq 2l$ since $\wt (a_{-2l-1}b)+2l=(\wt a+2l)+(\wt b+2l)=p+q$ if $\wt a+2l=p$ and $\wt b+2l=q$. 

Finally, by the skew-symmetry of the vertex operator, we have 
	\begin{align*}
	a_{-2l-1}b&=\Res_z Y(a,z)bz^{-2l-1}=\Res_z e^{zL(-1)}Y(b,-z)a z^{-2l-1}\\
	&=\Res_z \sum_{j\geq 0} \frac{L(-1)^j}{j!} z^j \sum_{n\in \Z} (-1)^{n+1} z^{-n-1-2l-1} b_n a\\
	&=b_{-2l-1}a+\sum_{j\geq 1} \frac{L(-1)^j}{j!} (-1)^j  b_{j-2l-1}a\\
	&\equiv b_{-2l-1}a+\sum_{j=1}^{2l} (-1)^j (b_{j-2l-1}a)_{-1-j}\vac \pmod{C_{2l+2}(V)}.
	\end{align*}
Thus, $R_{2l+2}(V)$ is commutative if and only if the obstruction term $\sum_{j=1}^{2l} (-1)^j (b_{j-2l-1}a)_{-1-j}\vac$ is in $C_{2l+2}(V)$. 
	\end{proof}
\begin{remark}
	If $l=0$, then $R_{2l+2}(V)=R_2(V)$ in view of \eqref{prodR2l+2V} and \eqref{RVprod}. The obstruction term for the commutativity in Lemma~\ref{lmR2l+2V} does not exist in this case. Hence $R_2(V)$ is commutative. 
	\end{remark}

We wish to find a $g$-twisted higher order analog of the epimorphism $\varphi$ in \eqref{varphi}. However, it turns out that $\varphi$ is not always generalizable without any extra assumptions.

\begin{thm}\label{mainthm2}
Let $V$ be a VOA, $g\in \mathrm{Aut}(V)$ be of order $T$, and $n=l+\frac{i}{T}\in \frac{1}{T}\Z$, where $l\in \N$ and $0\leq i\leq T-1$. Then there is a surjective linear map:
\begin{equation}\label{varphin}
\begin{aligned}
\varphi_n: R_{2l+2}(V)&\ra \gr A_{g,n}(V)=\bigoplus_{p=2l}^{\infty}F_pA_{g,n}(V)/F_{p-1}A_{g,n}(V),\\
a+C_{2l+2}(V)&\mapsto [a]+F_{p-1}A_{g,n}(V),\quad a+C_{2l+2}(V)\in R_{2l+2}(V)_p.
\end{aligned}
\end{equation}
If, furthermore, $i< \lfloor T/2\rfloor$, then $\varphi_n$ is an epimorphism of associative algebras. 
	\end{thm}
\begin{proof}
	Similar to Theorem~\ref{thm:AgVN}, we first show $\varphi_n$ is well-defined. Let $a\in V_{p-2l}\cap V^r$ with $0\leq r\leq T-1$ and $b\in V_{q-2l}$, for some $p,q\geq 2l$. Then $a_{-2l-2}b\in V_{p+q+1-2l}$, and $\varphi_n(a_{-2l-2}b+C_{2l+2}(V))=[a_{-2l-2}b]+F_{p+q}A_{g,n}(V)$. 
	 We need to show $[a_{-2l-2}b]\equiv 0\pmod{F_{p+q}A_{g,n}(V)}$. 
	 
	 Indeed, for  by Lemma 2.2 in \cite{DLM98(2)}, we have 
	 \begin{equation}\label{useful2}
	 \Res_z Y(a,z)b\frac{(1+z)^{\wt a-1+\delta_{i}(r)+l+r/T}}{z^{2l+\delta_i(r)+\delta_i(T-r)+m}}\in O_{g,n}(V),
	 \end{equation}
	 for any $m\geq 0$. Since $\delta_i(r)$ and $\delta_i(T-r)$ are either $0$ or $1$ in view of \eqref{AgnVprod}, we may choose $m\geq 0$ in such a way that $2l+\delta_i(r)+\delta_i(T-r)+m=2l+2$. Then by \eqref{useful2}, 
	 $$[a_{-2l-2}b]=-\sum_{j\geq 1} \binom{\wt a-1+\delta_i(r)+l+r/T}{j} [a_{j-2l-2}b]\in F_{p+q}A_{g,n}(V)$$
	 since $\wt (a_{j-2l-2}b)+2l=(p-2l)-j+2l+1+(q-2l)+2l\leq p+q$ for any $j\geq 1$, which means $[a_{j-2l-2}b]\in F_{p+q}A_{g,n}(V)$ by \eqref{levelfilonAgnV}. This proves the well-definedness of $\varphi_n$.
	 
	 Clearly, $\varphi_n$ is surjective. For any $a+C_{2l+2}(V)\in R_{2l+2}(V)_{p}$, with $\wt a+2l=p$, we have $[a]\in F_{p}A_{g,n}(V)$ by \eqref{levelfilonAgnV}. Hence $\varphi_{n}(R_{2l+2}(V)_{p})\ssq F_pA_{g,n}(V)/F_{p-1}A_{g,n}(V)$ for any $p\geq 2l$. 
	
	Finally, we show $\varphi_n$ is a homomorphism when $i< \lfloor T/2\rfloor$. Again, we let $a\in V_{p-2l}\cap V^r$ with $0\leq r\leq T-1$ and $b\in V_{q-2l}$. If $r=0$, since $\wt (a_{-2l-1}b)+2l=p+q$, then by \eqref{prodgrAgnV} and \eqref{prodR2l+2V}, 
	\begin{align*}
	&\varphi_n ((a+C_{2l+2}(V))\cdot (b+C_{2l+2}(V)))=\varphi_n ((-1)^l\binom{2l}{l}a_{-2l-1}b+C_{2l+2}(V))\\
	&=(-1)^l \binom{2l}{l} [a_{-2l-1}b]+F_{p+q-1}A_{g,n}(V)= \left([a]+F_{p-1}A_{g,n}(V)\right)\ast _{g,n} \left([b]+F_{q-1}A_{g,n}(V)\right)\\
	&= \varphi_n(a+C_{2l+2}(V))\ast_{g,n} \varphi_n(b+C_{2l+2}(V)).
	\end{align*}
	Now consider the case when $r>0$. Since $i< \lfloor T/2\rfloor$, the inequalities $r\leq i$ and $T-r\leq i$ cannot be satisfied simultaneously. Thus $\delta_i(r)+\delta_{i}(T-r)=1$ for any $r>0$. By \eqref{def:OgnV}, we have 
	$$[a_{-2l-1}b]=-\sum_{j\geq 1}\binom{\wt a-1+\delta_i(r)+l+r/T}{j} [a_{j-2l-1}b]\in F_{p+q-1}A_{g,n}(V),$$
	since $\wt (a_{j-2l-1}b)+2l=(p-2l)-j+2l+(q-2l)+2l\leq p+q-1$ for any $j\geq 1$. Then 
	\begin{align*}
		&\varphi_n ((a+C_{2l+2}(V))\cdot (b+C_{2l+2}(V)))	=(-1)^l \binom{2l}{l} [a_{-2l-1}b]+F_{p+q-1}A_{g,n}(V)\\
		&=0= \varphi_n(a+C_{2l+2}(V))\ast_{g,n} \varphi_n(b+C_{2l+2}(V)),
	\end{align*}
	in view of \eqref{prodgrAgnV}. 
	\end{proof}

In the proof of the modular invariance of $C_2$-cofinite VOAs, Miyamoto proved that $A_n(V)$ are finite-dimensional for all $n\geq 0$ if $V$ is $C_2$-cofinite, see Theorem 2.5 in \cite{M04}. The  following Corollary of Theorem~\ref{mainthm2} generalizes Miyamoto's result to the $g$-twisted case. 

\begin{coro}\label{coroofmain2}
Let $n=l+\frac{i}{T}\in \frac{1}{T}\Z$, where $l\in \N$ and $0\leq i\leq T-1$.
\begin{enumerate}
	\item If $V$ is $C_2$-cofinite, then $ A_{g,n}(V)$ is a finite-dimensional associative algebra.
	\item If $i< \lfloor T/2\rfloor$, and $R_{2l+2}(V)$ is a finitely generated associative algebra with respect to the product \eqref{prodR2l+2V}, then $ A_{g,n}(V)$ is left and right Noetherian.
	\end{enumerate}
	\end{coro}
\begin{proof}
	It was proved by Gaberdiel  and Neitzke that $V$ is $C_u$-cofinite for any $u\geq 2$ if $V$ is $C_2$-cofinite, see Theorem 11 in \cite{GN03}. Since the filtration \eqref{filAgnV} is exhaustive, by \eqref{varphin},
	$$\dim A_{g,n}(V)=\dim \gr A_{g,n}(V) \leq \dim R_{2l+2}(V)<\infty$$
	if $V$ is $C_2$-cofinite. Now let $i< \lfloor T/2\rfloor$, and assume that $R_{2l+2}(V)$ is a finitely generated algebra. By Theorem~\ref{mainthm2}, $\varphi_n$ is an epimorphism of associative algebras, then $\gr A_{g,n}(V)$ is a finitely generated commutative algebra in view of Remark~\ref{commofgrAgnV}, which is necessarily Noetherian. Hence $A_{g,n}(V)$ is left and right Noetherian by Proposition~\ref{prop2.2}. 
	\end{proof}

\begin{example}\label{ex:heisen}
Let $V=M_{\hat{\h}}(1,0)$ be the rank-one Heisenberg VOA, and let $g=\Id_V$. Then $V$ is $C_1$-cofinite by Theorem~\ref{equiC1}, since $R_2(V)\cong \C[x]$ as a commutative algebra, see \cite{DLM02}. In \cite{AB23}, Addabbo and Barron conjectured that for any $n\geq 1$, one has
\begin{equation}\label{AnVHeisenberg}
	A_n(V)\cong \mathrm{Mat}_{p(n)}(\C[x])\op A_{n-1}(V)
	\end{equation}
as a direct product of associative algebras, where $p(n)$ is the number of partitions of $n$. The isomorphism \eqref{AnVHeisenberg} was proved recently by Damiolini, Gibney, and Krashen, see Corollary 7.3.1 \cite{DGK23}. In particular, since $ \mathrm{Mat}_{p(n)}(\C[x])$ is a finitely generated module over a Noetherian ring $\C[x]$, using induction on $n$, it is easy to show that $A_n(V)$ are Noetherian for all $n\geq 0$. 
	\end{example}

We believe the following statement that generalizes Theorem~\ref{thm:AgVN} is true: 

\begin{conj}
Let $V$ be a VOA, $g\in \mathrm{Aut}(V)$ be of order $T$, and $n=l+\frac{i}{T}\in \frac{1}{T}\Z$, where $l\in \N$ and $i< \lfloor T/2\rfloor$. Then $A_{g,n}(V)$ is left and right Noetherian if $V$ is $C_1$-cofinite. 
 	\end{conj}

\subsection{Finiteness of $A_{g,n}(M)$} 
Buhl extended Gaberdiel and Neitzke's theorem to the case of $V$-modules. He proved that $A_{n}(M)$ are finite-dimensional for all $n\geq 0$ if $V$ is $C_2$-cofinite and $M$ is $C_2$-cofinite, see \cite{Bu02} Corollary 5.5. As another Corollary of Theorem~\ref{mainthm2}, we generalize Buhl's theorem to the twisted case. 
%Recall the following notations from \cite{Bu02}. Assume $V=X+C_2(V)$, where $\dim X<\infty$. i.e., $V$ is $C_2$-cofinite. Let $N>0$ be an integer such that $\bigoplus_{i>N} V_i\ssq C_2(V)$,  $B:=\max\{\wt x: x\in X\}$, and $Q:=\max\{N, 2B-1\}+1$. Let $M$ be a weak $V$-module. Given $w\in W$, let $L>0$ be the smallest positive integer such that $x_m w=0$ for all $x\in X$ and $m\geq L$. Note that $Q$ and $L$ are fixed numbers determined by $V$ and $w\in M$. 
 
% \begin{thm}[\cite{Bu02}]
%With the notations as above, let $M$ be a weak $V$-module generated by $w$. Then $M$ is spanned by elements of the following form: 
% \begin{equation}
% x^1_{-n_1}x^2_{-n_2}\ds x^k_{-n_k}w,\quad \mathrm{where}\quad x^i\in X,\ n_{1}\geq n_2\geq \ds \geq n_k\geq -L,
% \end{equation}
% where $n_1,\ds,n_k$ satisfy the following condition: If $n_j>0$ then $n_j>n_{j'}$ for $j<j'$; if $n_j\leq 0$ then $n_j=n_j'$ for at most $Q-1$ indices $j'$. 
% 	\end{thm}

 Abe, Buhl, and Dong proved that if $V$ is $C_2$-cofinite then an irreducible $V$-module $M$ is $C_2$-cofnite, see Proposition 5.2 in \cite{ABD03}. In fact, it is easy to show that $M$ is also $C_1$-cofinite by adopting a similar proof. Moreover, Buhl proved that $M$ is $C_n$-cofinite for all $n\geq 2$, if $M$ is $C_2$-cofinite, see Corollary 5.3 in \cite{Bu02}. Hence we have the following: 
 
 \begin{lm}\label{lastlm}
 	Let $V$ be $C_2$-cofinite, and let $M$ be an irreducible admissible $V$-module. Then $M$ is $C_n$-cofinite, for any $n\geq 1$.
 	\end{lm}

%\begin{lm}
% Let $V$ be $C_2$-cofinite, and let $M$ be an irreducible admissible $V$-module, then $M$ is $C_{n}$-cofinite for all $n\geq 2$. 
%\end{lm}
 
% He proved that if $V$ is $C_2$-cofinite, and an admissible $V$-module $M$ is $C_2$-cofinite, then $M$ is $C_n$-cofinite for all $n\geq 2$, see  Using this fact, Buhl also proved that $A_n(M)$ are finite-dimensional for all $n\geq 2$ under these conditions. 
 
We can generalize $\widetilde{C}^g_1(M)$ \eqref{weakC1gM} to the higher order case. For $l\geq 0$, define 
 \begin{equation}\label{generalC1}
 \widetilde{C}^g_{2l+1}(M):=\spn\{ a_{-2l-1}v:a\in V_+, v\in M\}+\spn\{b_0u: b\in \oplus_{p\geq 2}V_p\cap V^0, u\in M\}. 
 \end{equation}
 We say that $M$ is {\em weakly $C^g_{2l+1}$-cofinite} if $\dim M/\widetilde{C}^g_{2l+1}(M)<\infty$.

 \begin{coro}\label{coro2ofmain2}
 	Let $n=l+\frac{i}{T}\in \frac{1}{T}\Z$, where $l\in \N$ and $0\leq i\leq T-1$. 
 	
 	\begin{enumerate}
 		\item If $M=Y+\widetilde{C}^g_{2l+1}(M)$, $Y=\spn\{y^p\in:p\in \Lambda\}$, and $i< \lfloor T/2\rfloor$, then 
 		\begin{equation}\label{AgnMbigenerate}
 		A_{g,n}(M)=\sum_{p\in \Lambda} A_{g,n}(V)\ast_{g,n} [y^p]\ast_{g,n} A_{g,n}(V).
 		\end{equation}
 		In particular, if $M$ is weakly $C_{2l+1}^g$-cofinite, then $A_{g,n}(M)$ is a finitely generated $A_{g,n}(V)$-bimodule.
 \item If $M=U+C_{2l+1}(M)$, $U=\spn\{u^\al:\al\in I\}$, and $i< \lfloor T/2\rfloor$, then 
 	\begin{equation}\label{AgnMgenerating}
 	A_{g,n}(M)=\sum_{\al\in I} A_{g,n}(V)\ast _{g,n}[u^\al]=\sum_{\al\in I} [u^\al]\ast_{g,n} A_{g,n}(V). 
 	\end{equation}
 	In particular, if  if $M$ is $C_{2l+1}$-cofinite, then $A_{g,n}(M)$ is finitely generated as a left or right $A_{g,n}(V)$-module.
 	\item If $M=W+C_{2l+2}(M)$ and $W=\spn\{w^j:j\in J\}$, then 
 	\begin{equation}\label{AgnMspan}
 	A_{g,n}(M)=\sum_{j\in J} \C [w^j]. 
 	\end{equation}
 		In particular, if $V$ is $C_2$-cofinite, then $A_{g,n}(M)$ is finite-dimensional for all $n\geq 0$.
 	\end{enumerate}
%  If $M$ is $C_{2l+1}$-cofinite, then $A_{g,n}(M)$ is a finitely generated left or right $A_{g,n}(V)$-module. 
\end{coro}
\begin{proof}
The proof is similar to the proof of Theorem~\ref{thm:AgMfg}. We write out the details for \eqref{AgnMbigenerate} and omit the rests. Denote $\sum_{p\in \Lambda} A_{g,n}(V)\ast_{g,n}[y^p]\ast_{g,n}A_{g,n}(V)$ by $N$. We use induction on the degree $m$ of $M=\bigoplus_{m=0}^\infty M(m)$ to show $[M(m)]\ssq N$. Since $\deg (a_{-2l-1}v)>0$ and $\deg (b_0u)>0$ for any  $a\in V_+$, $b\in \oplus_{p\in 2}V_p$, and $u,v\in M$, we have $\widetilde{C}^g_{2l+1}(M)\cap M(0)=0$ in view of \eqref{generalC1}. Hence $[M(0)]\ssq N$. Suppose the conclusion holds for smaller $m$. For $x\in M(m)$, we may assume 
	$$x=u+\sum_{k=1}^s a^k_{-2l-1}v^k+\sum_{q=1}^t b^q_0 u^q ,\quad u\in U,\ a^k\in V_+\cap V^r,\ b\in \oplus_{p\geq 2}V_p\cap V^r,\ v^k,u^q\in M,$$
	where $0\leq r\leq T-1$, $\wt a^k+2l+\deg v^k=m$, and $\wt b^q-1+\deg u^q=m$ for all $k,q$. We need to show that $[a^k_{-2l-1}v^k]\in N$ and $[b^q_0 u^q]\in N$ for all $k,q$.
	
	Recall the following formulas for the definition of  $A_{g,n}(M)=M/O_{g,n}(M)$ in \cite{JJ16}: 
\begin{equation}\label{4.15}
a\circ_{g,n} v=\Res_z Y_M(a,z)v\frac{(1+z)^{\wt a+l-1+\delta_i(r)+\frac{r}{T}}}{z^{2l+\delta_i(r)+\delta_i(T-r)}},
\end{equation}
where $\delta_i$ is defined by \eqref{def:OgnV}. 
For $a\in V^r$ and $v\in M$, one has
\begin{equation}\label{4.16}
a\ast_{g,n}v=\Res_z \sum_{m=0}^l (-1)^m\binom{m+l}{l} \Res_z Y(a,z)v \frac{(1+z)^{\wt a+l}}{z^{l+m+1}}
\end{equation}
if $a\in V^0$, and $a\ast_{g,n}v=0$ if $a\in V^r$ with $r>0$. 

For each $1\leq k\leq s$, if $a^k\in V^0$, by \eqref{4.16} and induction hypothesis, we have
\begin{align*}
[a^k_{-2k-1}v^k]=-\sum_{j,m\geq 0, j-m>-l}(-1)^m \binom{m+l}{l} \binom{\wt a^k+l}{j} [a^k_{j-l-m-1} v^k]\in N,
\end{align*} 
since $\deg (a^k_{j-l-m-2} v^k)<\wt a^k+2l+\deg v^k=m$ when $j-m>-l$. On the other hand, if $a^k\in V^r$ with $r>0$, we have  $\delta_i(r)+\delta_i(T-r)=1$ since $i< \lfloor T/2\rfloor$. By \eqref{4.15} we have 
$$[0]=[a^k\circ_{g,n}v^k]=[a^k_{-2l-1}v^k]+\sum_{j\geq 1}\binom{\wt a^k+l-1+\delta_i(r)+r/T}{j} [a^k_{j-2l-1}v^k]$$
in $A_{g,n}(M)$. Since $\deg (a^k_{j-2l-1}v^k)=\wt a^k-j+2l+\deg v^k<m$, then $[a^k_{j-2l-1}v^k]\in N$ for all $j\geq 1$ by the induction hypothesis. Hence $[a^k_{-2l-1}v^k]\in N$. Moreover, we have  
\begin{equation}\label{4.20}
b\ast_{g,n}u-u\ast_{g,n}b\equiv \Res_{z} Y(b,z)u(1+z)^{\wt b-1}\pmod{O_{g,n}(M)}
\end{equation}
for any $b\in V^0$ and $u\in M$, see Lemma 3.1 \cite{JJ16}. Then by \eqref{4.20} and \eqref{3.6}, we have 
$$[b^q_0 u^q]=-\sum_{j\geq 1} \binom{\wt b^q-1}{j} [b^q_j u^q]+[b^q]\ast_{g,n} [u^q]-[u^q]\ast_{g,n} [b^q]\in N,$$
since $\deg (b^q_j u^q)=\wt b^q-j-1+\deg u^q<m$ for any $j\geq 1$. 
Thus $[x]=[u]+\sum_{k=1}^s [a^k_{-2l-1}v^k]+\sum_{q=1}^t[b^q_0u^q]\in N$.
This shows \eqref{AgnMbigenerate}. 

Finally, assume $V$ is $C_2$-cofinite. By Lemma~\ref{lastlm}, $M$ is $C_{2l+2}$-cofinite, for any $l\geq 0$. Then by \eqref{AgnMspan}, $A_{g,n}(M)$ is finite-dimensional. 
\end{proof}

\section{Acknowledgment}
I'm deeply grateful to professors Angela Gibney and Danny Karshen for their encouragement and many valuable comments and suggestions for this paper.

%----------------------------------------------------------------------------------------------------

%\section{Declarations}
%\subsection{Ethical Approval } Not applicable
%\subsection{Competing interests } There are no competing interests in our work.
%\subsection{Authors' contributions } This is the lead author's solo work.
%\subsection{Funding} We don't receive any funding for this work.
%\subsection{Availability of data and materials }Not applicable

\end{document}